\documentclass[preprint,12pt]{elsarticle}

\usepackage{epsfig}
\usepackage{float}
\usepackage{amsmath,amssymb}
\usepackage{amsthm}
\usepackage{algorithm,algorithmic}
\usepackage[labelformat=simple]{subfig}

\usepackage{latexsym,bm}        

\newtheorem{theorem}{Theorem}[section]

\newtheorem{remark}[theorem]{Remark}

\newcommand{\set}[2]{\left\{{#1}\,:~{#2}\right\}}
\newcommand {\average}[1] {\mbox{$\left\{\!\!\left\{ #1 \right\}\!\!\right\}$}}
\newcommand {\jump}[1] {\mbox{$\left[\!\left[ #1 \right]\!\right]$}}

\graphicspath{{./eps/}}

\makeatletter
\def\input@path{{SourceFiles/}}
\makeatother
\usepackage{hyperref}

\begin{document}

\begin{frontmatter}

\title{Structure preserving reduced order modeling for gradient systems}

\author[thk]{Tu\u{g}ba Akman Y{\i}ld{\i}z\corref{cor1}}
\ead{takman@thk.edu.tr,tr.tugba.akman@gmail.com}

\author[sin]{Murat Uzunca}
\ead{muzunca@sinop.edu.tr}

\author[iam]{B\"ulent Karas\"{o}zen}
\ead{bulent@metu.edu.tr}

\cortext[cor1]{Corresponding author}

\address[thk]{Department of Logistics Management, University of Turkish Aeronautical Association, 06790, Ankara, Turkey }
\address[sin]{Department of Mathematics, Sinop University, 57000, Sinop, Turkey}
\address[iam]{Institute of Applied Mathematics \& Department of Mathematics, Middle East Technical University, 06800 Ankara, Turkey}

\begin{abstract}

Minimization of energy in gradient systems leads to formation of oscillatory and Turing patterns in reaction-diffusion systems. These patterns should be accurately computed  using fine space and time meshes over long time horizons to reach the spatially inhomogeneous   steady state. In this paper, a reduced order model (ROM) is developed  which preserves the gradient dissipative  structure. The coupled system of reaction-diffusion equations are discretized in space by the symmetric interior penalty discontinuous Galerkin (SIPG) method. The resulting system of ordinary differential equations (ODEs)  are integrated in time by the average vector field (AVF) method, which preserves the energy dissipation of the gradient systems. The ROMs are  constructed by the proper orthogonal decomposition (POD) with Galerkin projection.  The nonlinear reaction terms are computed  efficiently by discrete empirical interpolation method (DEIM). Preservation of the discrete energy  of the FOMs and ROMs with POD-DEIM  ensures the long term stability of the steady state solutions.
Numerical simulations are performed for the gradient dissipative systems with two specific equations; real Ginzburg-Landau equation and Swift-Hohenberg equation. Numerical results demonstrate that the POD-DEIM reduced order solutions preserve well the energy dissipation over time and at the steady state.

\end{abstract}


\end{frontmatter}

\section{Introduction}

It is well known that energy minimization in gradient systems leads to formation of patterns. The energy driven pattern formation occurs in form of  oscillatory and Turing patterns in gradient systems \cite{Kohn07,Kuwamura03, Kuwamura05,Kuwamura07}. The long-term accurate computation of  spatially heterogeneous steady state  patterns require fine space and time meshes. Due to the implicit nature of time integrators used, a coupled fully nonlinear system has to be solved by the Newton-Raphson method at each time step with high accuracy. The discrete energy  should be computed with high accuracy, when the nonlinear system is solved  up to the machine precision. Therefore, the numerical solution of gradient systems for accurate resolution of patterns over long time integration might be very time consuming, especially for two and three dimensional (2D and 3D) problems.

The goal of this work is to develop an accurate and efficient reduced order model of 2D gradient systems by preserving the energy dissipation. During the last decade, model order reduction techniques  have emerged as a powerful tool to reduce the cost of evaluating large systems of partial differential equations (PDEs). A low-dimensional linear reduced space is constructed that approximately represents the solution to the system of PDEs. Then, the solution can be obtained with a significantly reduced computational complexity  by projection of the original system onto
the reduced space (by Galerkin projection). Model reduction methods should
produce low-dimensional surrogate models capable of mimicking the structural features of the underlying physics
of the original system model. Reduced models that do not inherit key structural features and conservation properties of the original system
may produce unphysical solutions.

The proper orthogonal decomposition (POD)  is the most commonly used reduced order modeling technique in large-scale numerical simulations of complex systems. The stability of reduced order models over long-time integration and the structure preserving properties have been recently investigated in the context of Lagrangian systems \cite{Lall03,Carlberg15}, and for port-Hamiltonian systems \cite{Beattie11,Chaturantabu16}. For Hamiltonian and dissipative Hamiltonian systems,  symplectic model reduction techniques are constructed in \cite{Afkham18,Hesthaven16,Peng16}.
The average vector field (AVF) method  is applied as energy preserving/dissipating  time integrator to Hamiltonian and dissipative PDEs \cite{Celledoni12ped} like the nonlinear Schr{\"o}dinger equation, Korteweg-de  Vries equation, Sine-Gordon equation, Allen-Cahn equation, Ginzburg Landau equation and Maxwell equation. Applying AVF  with reduced order modeling leads to stable and efficient solutions for  energy dissipating systems like Allen-Cahn equation \cite{Karasozen17}, FitzHugh-Nagumo equation \cite{Karasozen17a} and energy conserving Hamiltonian systems \cite{Gong14,Karasozen18}. There are a few works in the literature about the reduced order modeling of pattern formation in reaction-diffusion systems, like the complex Ginzburg-Landau equation \cite{Ilak10, Terregani14} and FitzHugh Nagumo equation \cite{Karasozen17a}.

For the space discretization, we apply symmetric interior penalty discontinuous Galerkin (SIPG) method \cite{Arnold02,riviere08dgm}, and for time discretization, we utilize the energy preserving AVF method. The full order model based on SIPG-AVF discretization in space and time preserves the gradient structure of the Allen-Cahn equation \cite{Karasoen18a}, Cahn-Hilliard equation \cite{Karasoen17b} and the skew-gradient structure of the FitzHugh-Nagumo equation \cite{Karasozen17a}. The POD-Galerkin approach  involves an offline-online splitting strategy. In the offline stage, the high fidelity or truth solutions are  generated by numerical simulations of the discretized high dimensional full order model (FOM). The POD is then  applied to compute an optimal subspace to fit the high fidelity data. A reduced system is constructed by projecting the high-dimensional system to this subspace. In the online stage, the reduced system is solved in the low-dimensional subspace. The primary challenge in producing the low dimensional models of the high dimensional discretized PDEs is the efficient evaluation of the nonlinearities on the POD basis. 
 Here, we apply the  discrete empirical interpolation (DEIM)  \cite{chaturantabut10nmr} with the SIPG space discretization which is more efficient for evaluations of the nonlinear terms than the continuous FEM \cite{Karasozen17a} in reduced order models. 
The reduced order basis functions in POD-DEIM are computed in the offline stage using  the singular value decomposition (SVD), which can be computationally demanding for large snapshot matrices. Here, we use the randomized singular value decomposition (rSVD) \cite{Halko11a,Mahoney11,Martinsson16} as  a fast and accurate alternative to the deterministic SVD  to reduce the computational cost in the offline stage \cite{Erichson16}. 

In the literature, there are some recent works concerning the ROMs based on POD which preserve energy dissipative structure of the model. For example, the study \cite{Karasozen18} offers mass and energy preserving FOM and ROM for two dimensional nonlinear Schr{\"o}dinger equation (NLSE). To tackle the nonlinearity, DEIM is applied and it is shown that the discrete energy is preserved approximately for DEIM-ROM approximation. NLSE is solved by decomposing the real and complex parts of the solution which leads to a system of nonlinear PDEs. However, it is a Hamiltonian ODE. In another work \cite{Karasozen17}, an energy stable ROM for a parametrized Allen-Cahn equation is proposed based on POD-greedy adaptive sampling of the snapshots in time and the greedy DEIM for the nonlinear term. It is reported that the energy decreases almost unconditionally for POD-DEIM approximation with large time step-size. On the other hand, a skew-gradient system, namely FitzHugh-Nagumo equation, is investigated and a ROM-DEIM approximation is obtained and compared with FOM \cite{Karasozen17a}. It is justified that the mini-maximizer property of the continuous energy holds for the full-discrete model. As different from the works mentioned above, in this current study, we investigate the structure preserving reduced order modeling for two typical gradient systems; Ginzburg-Landau equation and Swift-Hohenberg equation. Decomposing the solution of Ginzburg-Landau equation in real and complex parts, we obtain a system of parabolic-parabolic equations with cubic nonlinearities. On the other hand, the fourth order Swift-Hohenberg equation can be reduced to a system of parabolic-elliptic equations and the parabolic equation contains a cubic nonlinear term. After obtaining the systems, the full-discrete model is obtained by applying SIPG in space and AVF in time. Then, we derive the corresponding parabolic-parabolic and parabolic-elliptic ROMs with the use of DEIM for the nonlinearities. We prove the energy dissipation of the POD and POD-DEIM reduced order models and the contribution of the parabolic-parabolic and parabolic-elliptic forms is reflected to the energy analysis. In detail, energy decreasing property of the gradient systems are proven under some condition on time step-size and this condition is related to both of the solutions of the parabolic-parabolic system, while the bound on time step-size depends only on the solution of the parabolic PDE for the parabolic-elliptic system. Moreover, we mention a relation between FOM and ROM so that the energy dissipation of FOM can be deduced easily. Furthermore, we compare DEIM and POD with respect to accuracy and efficiency; the POD-DEIM solutions have almost the same accuracy as the POD, but POD-DEIM is faster than POD. The steady state solutions of the gradient possesses Hamiltonian structure \cite{Kuwamura03, Kuwamura05,Kuwamura07}. Our numerical results show that around  the steady state, the energy of the FOMs and reduced order models (ROMs) remain constant. We drive upper bounds for the time meshes for the POD-DEIM reduced order models of the two gradient systems above. To the best of our knowledge, this is the first study justifying energy dissipation of POD-DEIM approximation for two typical example of gradient systems; Ginzburg-Landau equation and Swift-Hohenberg equation.

The remainder of the paper is organized as follows. The gradient systems are reviewed briefly in Section~\ref{Sec:grad}. Space-time discretization by SIPG and AVF method is given in Section~\ref{Sec:fom}. In Section~\ref{Sec:rom}, a small dimensional ROM by utilizing the POD method and also the POD-DEIM are presented. In  Section~\ref{Sec:analysis},  we investigate the gradient dissipative structures of reduced order RGL and SH systems based on DEIM approximation. The preservation of the dissipative structure, efficiency and accuracy of the ROMs are illustrated through the numerical simulations in Section~\ref{Sec:numeric}. The paper ends with some conclusions.

\section{Gradient systems}
\label{Sec:grad}

In the following section, we summarize the properties of gradient systems \cite{Kuwamura03,Kuwamura05,Yanigada02mfrd}. Let us consider an $m$-component reaction-diffusion system
\begin{equation} \label{grad}
S \partial_t \vec{w} = D \Delta \vec{w} + {\bm f}(\vec{w}), \qquad
{\bm f}(\vec{w}) = Q \nabla _{\vec{w}} F(\vec{w}), \quad \hbox{in } \Omega \times [0,T],
\end{equation}
where $\vec{w}({\bm x},t) = (w_1({\bm x},t), \ldots, w_m({\bm x},t))^T$ and $\Delta$ is the Laplace operator acting component-wise on $\vec{w}$, $t$ is the time variable, $T > 0$ is the final time and $\Omega \in \mathbb{R}^d\; (d=1,2,3)$ is a bounded polygonal domain. We consider two-dimensional reaction diffusion equations ($d=2$) and set the spatial variable as ${\bm x} = (x,y)$. In \eqref{grad}, $F(\vec{w})$ is a real-valued potential functional, $S$ is a nonnegative diagonal matrix, and $Q$ is a non-degenerate symmetric matrix with $Q^2 =I$ and satisfies the condition $D^T Q = Q D$, which ensures the nondegeneracy of $QD$. For the nonlinear vector ${\bm f}(\vec{w})=(f_1(\vec{w}),\ldots , f_m(\vec{w}))^T$, the Jacobian matrix ${\bm f}_{\vec{w}}$ satisfies
$$
{\bm f}_{\vec{w}}(\vec{w})^TQ = Q{\bm f}_{\vec{w}}(\vec{w}).
$$
For the system \eqref{grad}, the corresponding energy functional is
\begin{equation} \label{grad_energy}
E[\vec{w}] = \int_{\Omega}  \left \{   \frac{1}{2} \langle D\nabla \vec{w}, Q \nabla \vec{w} \rangle_{\mathbb{R}^n} - F(\vec{w})   \right \} d{\bm x},
\end{equation}
with $\langle \cdot , \cdot \rangle_{\mathbb{R}^n} $ denoting the standard Euclidean inner product on $\mathbb{R}^n$, and it holds
$$
\frac{d}{dt} E[\vec{w}({\bm x},t)] = - \int_{\Omega} \langle\vec{w}_t,QS\vec{w}_t\rangle_{\mathbb{R}^n} d{\bm x} .
$$
The system \eqref{grad} has a  gradient structure when $QS$ is a nonnegative symmetric matrix. In gradient systems, like the real Ginzburg-Landau equation and Swift-Hohenberg equation, \eqref{grad_energy} is the free energy  which decreases monotonically,  and local minimizers of the free energy correspond to the stable steady states.

Below, we give two representative examples of the gradient systems, which are investigated in this paper. \\

\noindent {\em Real Ginzburg-Landau (RGL) equation  or  Ginzburg-Pitaevski equation\/}
\cite{Cross09,Hecke94,Hoyle06,Kuwamura03,Kuwamura05,levermore96,Saarlos95,Walgraef97} is given by
\begin{equation*}
\partial_t \psi = \Delta \psi  +  \mu \psi -|\psi|^2\psi,
\end{equation*}
where the unknown solution $\psi$ is complex-valued. Decomposing the solution $\psi = u + iv$ in real and complex parts yields the coupled system of reaction-diffusion equations
\begin{equation} \label{rgle}
\partial_t u   =  \Delta u  + \mu u  - (u^2 +v^2) u, \qquad
\partial_t v   =  \Delta v  + \mu v  - (u^2 +v^2) v,
\end{equation}
and it is of the form of the system \eqref{grad} through the setting $S=D=Q=I_2$, where $I_2$ is the two dimensional identity matrix, $\vec{w}=(u,v)^T$, and with the potential function $ F(u,v) = \mu (u^2 + v^2)/2 - (u^2+v^2)^2/4$.\\

\noindent {\em Swift-Hohenberg (SH) equation\/}  \cite{Deghan17,Gomez12,Kuwamura03,Kuwamura05,Moreno14,Peletier03,Peletier04,Perez14} is given by
\begin{equation*}
\partial_t u = \mu u - (1 +\Delta)^2u - u^3.
\end{equation*}
Putting $v = u + \Delta u$  yields the coupled system
\begin{equation}\label{she}
\partial_t u   =  -\Delta v    + \mu u -u^3 -v, \qquad
0   = \Delta u + u - v.
\end{equation}
In the sense of the system \eqref{grad}, the above system \eqref{she} is related to the setting
$$
S = \left ( \begin{array}{cc} 1 & 0\\ 0 & 0 \end{array} \right), \quad
D = \left ( \begin{array}{cc} 0 & -1\\ 1 & 0 \end{array} \right), \quad
Q = \left ( \begin{array}{cc} 1 & 0\\ 0 & -1\end{array} \right),
$$
with $\vec{w}=(u,v)^T$, and with the potential $ F(u,v) = 0.5\mu u^2 - 0.25u^4 -uv + 0.5v^2.$\\

The models \eqref{rgle} and \eqref{she} have cubic nonlinear reaction terms.

Turing patterns are formed, when the critical wave number near the stationary solution is non-zero. In this case $QD$ is indefinite. \cite{Kuwamura03,Kuwamura05}. For RGL equation, $QD$ is positive definite, so there exist no Turing patterns. However, there exist spatially periodic patterns around the stationary solutions \cite{Kuwamura03,Kuwamura05}. Turing instability conditions and selection of Turing patterns for gradient systems are investigated in \cite{Kuwamura05}.

Steady state solutions  of the gradient systems  admit Hamiltonian structure
\cite{Kuwamura03,Kuwamura05,Kuwamura07}. The stationary system
$$
0 = D \Delta \vec{w} + {\bm f}(\vec{w})
$$
can be written in Hamiltonian form as
$$
J \nabla \vec{z} =  \frac{\partial H(\vec{z})}{\partial \vec{z} }
$$
with the skew-symmetric matrix $J$ and the Hamiltonian given by
$$
J = \left ( \begin{array}{cc} 0 & -QD \\ QD & 0 \end{array} \right ), \qquad H(\vec{z}) := \frac{1}{2}\langle D\nabla \vec{w}, Q  \nabla \vec{w}\rangle + F(\vec{w}),
$$
where $\vec{z} = (\vec{w}, \nabla \vec{w} )^T$.

The RGL equation \eqref{rgle} and SH equation \eqref{she} have dissipative structures for homogeneous Dirichlet and periodic boundary conditions
\cite{Hecke94,levermore96,Matsuo01,Saarlos95}. Here, we consider periodic boundary conditions in both components for RGL and SH equations. 

\section{Discretization of the system (high-fidelity model)}
\label{Sec:fom}

In this section, we give a brief construction to the space-time discretization of the generic model \eqref{grad}, which we call high-fidelity or full order model (FOM). To do this, let us first consider the dimensionless scalar case of the generic model \eqref{grad} as
\begin{equation} \label{grad1}
\begin{aligned}
\partial_t w &= \Delta w + f(w), & \hbox{in } \Omega \times (0,T],\\
w(0,{\bm x}) &= w_0({\bm x}) , & f(w) = \nabla _{w} F(w),
\end{aligned}
\end{equation}
where $\Omega\subset\mathbb{R}^2$ is a bounded polygonal domain, and $w_0({\bm x})$ is a prescribed initial condition. Later, we extend to the generic model \eqref{grad} in a straightforward manner. In the sequel, we first give the semi-discretization of the system \eqref{grad1} by using SIPG method, which is a well-known member of the family of discontinuous Galerkin (dG) methods, for space discretization, and then we present the full discretization using AVF time integrator.

\subsection{Discontinuous Galerkin discretization in space}

Discontinuous Galerkin methods  gained an increasing importance for an efficient and accurate solution of PDEs in the last twenty years.  For the discretization the gradient systems in space we have the interior point discontinuous Galerkin method 
\cite{Arnold02, riviere08dgm}.

To obtain the infinite dimensional variational formulation, as usual, we multiply the equation \eqref{grad1} by a test function $\phi ({\bm x})\in W$, we take integral over the domain $\Omega$, and we apply Green's formula to the second order term. Then, the continuous weak formulation of the problem \eqref{grad1} can be written as follows: for a.e. $t\in(0,T]$, find $w:=w(t,\bm{x}) \in W$ such that for any $\phi:=\phi(\bm{x}) \in W$
\begin{equation}\label{discrete_grad}
	\begin{aligned}
	\langle \partial_t w,\phi\rangle_{L^2(\Omega)} &= -\langle\nabla w,\nabla\phi\rangle_{L^2(\Omega)}  +  \langle f(w),\phi\rangle_{L^2(\Omega)},\\
	\langle w,\phi \rangle_{L^2(\Omega)} &= \langle w_0,\phi \rangle_{L^2(\Omega)},
	\end{aligned}
\end{equation}
where $\langle\cdot , \cdot\rangle_{L^2(\Omega)}$ denotes the standard $L^2$-inner product on the space $L^2(\Omega)$ of square integrable functionals over the domain $\Omega\subset\mathbb{R}^2$, together with the induced $L^2$-norm $\|w\|_{L^2(\Omega)}=\sqrt{\langle w , w\rangle_{L^2(\Omega)}}$.
Throughout this paper, unless otherwise stated, $\langle \cdot , \cdot\rangle$ and $\|\cdot\|$ denote the standard $L^2$-inner product and $L^2$-norm over the aforementioned domains, respectively. The solution space $W$ is taken as the subset of the Sobolev space $H^{1}(\Omega)$ satisfying the periodicity on the boundary $\partial \Omega$ of the domain $\Omega$:
\begin{align*}
W &= \{\phi\in H^{1}(\Omega) | \; \phi \text{ periodic on }  \partial \Omega\},\\
H^{1}(\Omega) &= \{\phi\in L^{2}(\Omega) | \; \frac{\partial\phi}{\partial x}, \; \frac{\partial\phi}{\partial y}\in L^2(\Omega)\}.
\end{align*}

For the finite dimensional problem, continuous FEMs aim to approximate the solutions $w$ of \eqref{discrete_grad} from a conforming, finite dimensional subspace $W_h \subset W$. On the other hand, here we apply SIPG and we point out  that  the space of solutions/test functions consist of piecewise discontinuous polynomials in dG methods. That is, no continuity constraints are explicitly imposed on the state and test functions across the element interfaces. As a consequence, weak formulations must include jump terms across interfaces, and typically penalty terms are added to control the jump terms.

For discrete setting, we consider $\{ \mathcal{T}_h\}_h$ as a family of shape-regular simplicial triangulations of $\Omega$ for $h>0$. Each mesh $\mathcal{T}_h$ consists of closed triangles such that $\overline{\Omega} = \bigcup_{K \in \mathcal{T}_h} \overline{K}$ holds. The diameter of an (triangular) element $K$ and the length of an edge $E$ are denoted by $h_{K}$ and $h_E$, respectively. We split the set of all edges $\mathcal{E}_h$ into the set $\mathcal{E}^{0}_{h}$ of interior edges and the set $\mathcal{E}^{p}_{h}$ of (periodic) boundary edge-pairs. An individual element of the set $\mathcal{E}^{p}_{h}$ is of the form $\omega =\{E_l, E_m\}$ where $E_l \subset \partial K_{n_l}\cap \partial\Omega$, and $E_m \subset \partial K_{n_m}\cap \partial\Omega$ is the corresponding periodic edge-pair of $E_l$. 
We define the space of solution and test functions by
\begin{align}\label{dgspace}
W_h &=\set{w \in L^2(\Omega)}{ w\mid_{K}\in \mathbb{P}_q(K) \quad \forall K \in \mathcal{T}_h},
\end{align}
where $\mathbb{P}_q(K)$  is the set of polynomials  of degree at most $q$ in $K$, and that the space $W_h$ is a non-conforming space, i.e., $W_h \not\subset W$.
Let $w_h(0,\bm{x})\in W_h$ be the projection (orthogonal $L^2$-projection) of the initial condition $w_0(\bm{x})$ onto $W_h$.
Then, the SIPG weak formulation of the system \eqref{discrete_grad} reads as: for a.e. $t \in (0,T]$, find $w_h:=w_h(t,\bm{x}) \in W_h$ such that
\begin{align}\label{semi_dg}
	\langle \partial_t w_h, \phi \rangle  &=   -a_h(w_h, \phi) + b_{h}(w_h,\phi),  &  \forall \phi \in W_h.
\end{align}
In \eqref{semi_dg}, $a_h(\cdot , \cdot )$ and $b_{h}(\cdot , \cdot )$ are bi-linear and non-linear forms, respectively, which are given for any $\phi\in W_h$ by
\begin{align*}
	a_h(w_h,\phi)=& \sum \limits_{K \in \mathcal{T}_h} \int \limits_{K}  \big(  \nabla w_h \cdot  \nabla \phi )\; d\bm{x}
	+ \sum \limits_{ E \in \mathcal{E}^0_h} \frac{\kappa }{h_E} \int \limits_E \jump{w_h} \cdot \jump{\phi} \; d\bm{s}\\
	& -  \sum \limits_{ E \in \mathcal{E}^0_h } \int \limits_E \Big( \average{  \nabla w_h} \cdot \jump{\phi}   +  \average{  \nabla \phi} \cdot \jump{w_h} \Big) \; d\bm{s} + J^p_h(w_h,\phi), \\
	J^p_h(w_h,\phi) =& - \sum_{\omega \in \mathcal{E}^{p}_{h}} \int \limits_{\omega} \Big( \average{   \nabla w_h}_{\omega} \cdot \jump{\phi}_{\omega} + \average{   \nabla \phi}_{\omega} \cdot \jump{w_h}_{\omega} \Big) \; d\bm{s} \\
	& \quad  + \sum_{\omega \in \mathcal{E}^{p}_{h}} \frac{\kappa }{h_{E}} \int \limits_{\omega} \jump{w_h}_{\omega} \cdot \jump{\phi}_{\omega} \; d\bm{s}, \\
	b_{h}(w_h,\phi) =&  \sum \limits_{K \in \mathcal{T}_h} \int \limits_{K} f(w_h)\phi \; d\bm{x},
\end{align*}
where the bilinear form $J^p_h(w_h,\phi)$ is related to the periodic boundary edges.
Moreover, $\jump{\cdot}$ and $\average{\cdot}$ denote the jump and average operators in dG schemes, and the parameter $\kappa$ is called penalty parameter which should be sufficiently large to ensure the stability of the SIPG discretization with a lower bound depending only on the polynomial degree $q$, for details see \cite{riviere08dgm,Vemaganti07}.

The SIPG discretized semi-discrete solution of \eqref{semi_dg} is given by
\begin{equation}\label{dg_coef}
w_h(t,\bm{x}) = \sum \limits_{i=1}^{n_K} \sum \limits_{j=1}^{n_{q}} w_{j}^{i}(t) \varphi_{j}^{i}(\bm{x}),
\end{equation}
where $w_{j}^{i}(t)$ and $\varphi_{j}^{i}(\bm{x})$ are the unknown coefficients and the basis functions in $W_h$, respectively, for $j=1,2,\cdots, n_{q}$ and $i=1,2, \cdots, n_K$. The number $n_K$ denotes the number of (triangular) elements in $\mathcal{T}_h$, and $n_{q}$ is the local dimension on each element with the identity $n_{q}= (q+1)(q+2)/2$, where $q$ is the degree of the polynomial order. Note that the degrees of freedom in dG methods are given by $N:=n_K\times n_q$, and throughout this paper we denote the dimension of the high-fidelity model or FOM  by $N$. Inserting the expansions \eqref{dg_coef} into the system \eqref{semi_dg}, we obtain the following dynamical system:
\begin{equation}\label{fom}
\bm{M}\bm{w}_t = -\bm{A}\bm{w} + \bm{b}(\bm{w}),
\end{equation}
where ${\bm w}:={\bm w}(t) \in\mathbb{R}^N$ is the unknown coefficient vector for the solution $w_h$ with the ordered entries
\begin{equation*}
{\bm w} = (w_{1}^{1}(t), \cdots, w^{1}_{n_q}(t), \cdots, w_{n_K}^{1}(t), \cdots, w^{n_K}_{n_q}(t)),
\end{equation*}
the matrix $\bm{M}\in\mathbb{R}^{N\times N}$ is the usual mass matrix, $\bm{A}\in\mathbb{R}^{N\times N}$ is the stiffness matrix corresponding to the bilinear form $a_h(\cdot,\cdot)$, and $\bm{b}(\bm{w})\in\mathbb{R}^N$ is the vector corresponding to the non-linear form $b_{h}(\cdot,\cdot)$ with the entries $(\bm{b}(\bm{w}))_i = b_{h}(w_h,\varphi_i)$ for $i=1,\ldots , N$, where we have used the same ordering for basis functions as for the unknown coefficients:
\begin{equation*}
{\bm \varphi}(\bm{x}) = (\varphi_{1}^{1}(\bm{x}), \cdots, \varphi^{1}_{n_q}(\bm{x}), \cdots, \varphi_{n_K}^{1}(\bm{x}), \cdots, \varphi^{n_K}_{n_q}(\bm{x})).
\end{equation*}

\subsection{Time discretization}

Energy stable time discretization methods preserve the dissipative structure of the numerical solution of gradient flow equations.
Implicit Euler method and average vector field (AVF) method are energy stable time discretization techniques which are robust with small diffusion parameters. The AVF method is the only second order implicit energy
stable method and it preserves energy decreasing property for the gradient systems.

We split the time interval $[0, T]$ into $J$ equal length subintervals $(t_{k-1},t_k]$ as $0 = t_0 < t_1 < \ldots < t_J = T$ with the uniform step-size $\Delta t=t_k - t_{k-1},\; k= 1, 2, \ldots , J$.
We set $\bm{w}_{n} \approx \bm{w} (t_n)$ as the approximate solution vector at the time instance $t=t_n$. For $t=0$, let $w_h(0,\bm{x})\in W_h$ be the projection (orthogonal $L^2$-projection) of the initial condition $w_0({\bm x})$ onto $W_h$, and let ${\bm w}_0$ be the corresponding initial coefficient vector satisfying the expansion \eqref{dg_coef}. 
Then, applying the AVF method to the semi-discrete system \eqref{fom}, the full discrete problem of \eqref{grad1} is read as: for $n=0,1, \ldots , J-1$, find $\bm{w}_{n+1}$ satisfying
\begin{equation*}
\bm{M}\left(\frac{\bm{w}_{n+1}-\bm{w}_{n}}{\Delta t}\right)  = -\frac{1}{2}\bm{A}(\bm{w}_{n+1}+\bm{w}_{n}) +  \int _{0}^{1} \bm{b}(\xi \bm{w}_{n+1}+( 1-\xi )\bm{w}_{n}) \; d \xi.
\end{equation*}

\subsection{Full discretization of the model}

In this section, we extend the full discrete formulation of the dimensionless scalar model \eqref{grad1} to the generic model \eqref{grad}, and also we give the FOM formulations of the representative examples like RGL equation \eqref{rgle} and SH equation \eqref{she}.

Let us consider the generic model \eqref{grad} with $m$-component unknown solution $\vec{w}({\bm x},t)=(w_1({\bm x},t),\ldots , w_m({\bm x},t))^T$ together with its semi-discrete solution $\vec{w}_h({\bm x},t)=(w_{1,h}({\bm x},t),\ldots , w_{m,h}({\bm x},t))^T\in [W_h]^m$.  For each $i=1,\ldots , m$, we set $\bm{w}_{i,h}(t)\in\mathbb{R}^{N}$ as the coefficient vector of the semi-discrete solution $w_{i,h}({\bm x},t)$ satisfying the relation \eqref{dg_coef}, and we define the total coefficient vector as $\bm{w}(t)= (\bm{w}_{1,h}(t)^T,\ldots , \bm{w}_{m,h}(t)^T)^T\in\mathbb{R}^{mN}$. Moreover, for each $i=1,\ldots , m$, let $\bm{w}_{i,n}\approx \bm{w}_{i,h}(t_n)$ denote the approximate coefficient vectors at the time instance $t=t_n$, and approximate total coefficient vector is $\bm{w}_n \approx (\bm{w}_{1,n}^T,\ldots , \bm{w}_{m,n}^T)^T$.
 Finally, we set $\bm{w}_{0}=(\bm{w}_{1,0}^T,\ldots , \bm{w}_{m,0}^T)^T$, where each $\bm{w}_{i,0}$ is the $L^2$-projection of the initial conditions $w_{i,0}(\bm{x})$'s onto the dG space $W_h$. Then, the FOM of the model \eqref{grad} is given by
\begin{equation}\label{fom1}
(S\otimes\bm{M})\bm{w}_t = -(D\otimes\bm{A})\bm{w} + \bm{b}(\bm{w}),
\end{equation}
where $\otimes$ denotes the Kronecker product, and $\bm{b}(\bm{w})=(\bm{b}_1(\bm{w}), \ldots , \bm{b}_m(\bm{w}))^T\in\mathbb{R}^{mN}$ is the nonlinear vector whose entries are given by
$$
(\bm{b}_i(\bm{w}))_j=b_{i,h}(\vec{w}_h,\varphi_j) := \sum \limits_{K \in \mathcal{T}_h} \int \limits_{K} f_i(\vec{w}_h)\varphi_j \; d\bm{x}, \qquad j=1,\ldots , N.
$$

According to the general formulation \eqref{fom1}, and the related settings of $S$, $D$, $Q$ and potential functions $F(\vec{w})$ for the representative examples, their FOM formulations for the unknowns $\vec{w}:=(u,v)^T$, their coefficients vectors $\bm{w}:=(\bm{u},\bm{v})^T$, and the corresponding dG counter-part of the discrete energy functionals (the term $\langle D\nabla \vec{w}, Q \nabla \vec{w} \rangle_{\mathbb{R}^n}$ in the energy formulation \eqref{grad_energy} corresponds to the dG bilinear form $a_h(\cdot , \cdot)$ ) are expressed as\\

\noindent {\em FOM for RGL equation \eqref{rgle}\/}
\begin{equation}\label{fom_rgle}
\begin{aligned}
\bm{M}\bm{u}_t &= -\bm{A}\bm{u} + \bm{b}_1(\bm{u},\bm{v}), \quad f_1(u,v) =  \mu u - (u^2+v^2)u, \\
\bm{M}\bm{v}_t &= -\bm{A}\bm{v} + \bm{b}_2(\bm{u},\bm{v}), \quad f_2(u,v) =  \mu v - (u^2+v^2)v,\\
& E_h[u_h,v_h] = \frac{1}{2}a_h(u_h,u_h) + \frac{1}{2}a_h(v_h,v_h) - \langle F(u_h,v_h), 1 \rangle, \\
& F(u,v) = \frac{\mu}{2}(u^2+v^2) - \frac{1}{4}(u^2+v^2)^2.
\end{aligned}
\end{equation}

\noindent {\em FOM for SH equation \eqref{she}\/}
\begin{equation}\label{fom_she}
\begin{aligned}
\bm{M}\bm{u}_t &= \bm{A}\bm{v} + \bm{b}_1(\bm{u},\bm{v}), \quad f_1(u,v) =  \mu u - u^3 - v, \\
0 &= -\bm{A}\bm{u} + \bm{M}\bm{u} - \bm{M}\bm{v}, \quad f_2(u,v) =  u - v,\\
& E_h[u_h,v_h] = - a_h(u_h,v_h) - \langle F(u_h,v_h), 1 \rangle, \\
& F(u,v) = \frac{\mu}{2}u^2 - \frac{1}{4}u^4 - uv + \frac{1}{2}v^2.
\end{aligned}
\end{equation}

Note that in the SH equation \eqref{fom_she}, we do not require a nonlinear vector $\bm{b}_2(\bm{u},\bm{v})$ since the function $f_2(u,v)$ for the corresponding system is linear in both $u$ and $v$.

\section{Reduced order model}
\label{Sec:rom}

Because the computation of the FOM \eqref{fom1} is time consuming,  in this section, we construct a small dimensional ROM by utilizing the POD method \cite{Kunisch01}. We have constructed separate POD basis for $u$ and $v$  to preserve the structure of coupled system in reduced form. Similarly, the DEIM bases are constructed separately for the nonlinear vectors $\bm{b}_1$ and $\bm{b}_2$ in \eqref{fom_rgle}, and $\bm{b}_1$ in \eqref{fom_she}.
The low-rank approximation is computed in three steps: computation of the numerical solutions of the original high-fidelity system; dimensionality-reduction of the snapshot matrices by SVD;  Galerkin projection of the dynamics on the low-rank subspace. The first two steps are known as the offline stage, and the last one is the online stage. Offline stage is usually expensive and online step should be fast to run in real time.

As we did for dG formulation, we give, firstly, a brief construction procedure of the ROM for the scalar system \eqref{grad1} with the related dynamical system \eqref{fom}, and then we give the ROM formulations for our representative examples.

\subsection{POD Galerkin projection}

For the scalar system \eqref{grad1} with the $N$-dimensional FOM \eqref{fom}, the ROM of lower dimension $k\ll N$  is formed by the Galerkin projection of the system onto a $k$-dimensional reduced space
$$
W_h^r=\text{span} \{ \psi_{1}, \ldots, \psi_{k} \}\subset W_h,
$$
resulting in the reduced solution $w_h^r(t,\bm{x})\approx w_h(t,\bm{x})$, an approximation to the high-fidelity solution $w_h(t,\bm{x})$, and it  satisfies the weak formulation
\begin{align}\label{semi_dg_rom}
	\langle \partial_t w_h^r, \eta \rangle  &=   -a_h(w_h^r, \eta) + b_{h}(w_h^r,\eta),  &  \forall \eta \in W_h^r.
\end{align}

The reduced solution $w_h^r(t,\bm{x})$ is now of the form
\begin{equation}\label{rom_exp}
w_h^r(t,\bm{x}) = \sum_{i=1}^{k} w_i^r(t) \psi_{i}(\bm{x}),
\end{equation}
where $\psi_{i}(\bm{x})$'s are the orthogonal (in $L^2$-sense) reduced basis functions spanning the reduced space $W_h^r$, and  $w_i^r(t)$'s are the coefficients of the reduced solution, which we call reduced coefficients. From the reduced coefficients $w_i^r(t)$, we set the solution of the reduced system as $\bm{w}^r(t):= \left( w_1^r(t),\ldots, w_k^r(t) \right)^T$. Note that the reduced space $W_h^r\subset W_h$ is a subset of the dG space $W_h$; hence the reduced  basis functions are linear combination of the dG  basis functions $\{ \varphi_j \}_{j=1}^N$:
\begin{equation*}
\psi_{i}(\bm{x}) = \sum_{j=1}^{N} \Psi_{j,i} \varphi_j(\bm{x}), \qquad i=1,\ldots ,k,
\end{equation*}
where $\Psi_{j,i}$'s are the coefficients of the $i$th reduced basis function $\psi_{i}(\bm{x})$. Then, using the column vectors $\Psi_{i} = \left( \Psi_{1,i}, \ldots ,\Psi_{N,i} \right)^T$ of the coefficients, we construct the following matrix of POD modes
\begin{equation}
\bm{\Psi} := [ \Psi_{1}, \ldots, \Psi_{k}]\in \mathbb{R}^{N\times k}.
\end{equation}

To obtain the reduced basis functions $\{\psi_{i}(\bm{x}) \}_{i=1}^{k}$, we need to solve the minimization problem \cite{Kunisch01}
\begin{align*}
\min_{\psi_{1},\ldots ,\psi_{k}} \frac{1}{J}\sum_{j=1}^{J} \left\| w_h(t_j,\bm{x}) - \sum_{i=1}^k \langle w_h(t_j,\bm{x}),\psi_{i}(\bm{x})\rangle \psi_{i}(\bm{x})\right\|^2, \\
\text{subject to } \langle \psi_{i},\psi_{j}\rangle = \Psi_{i}^T\bm{M}\Psi_{j}=\delta_{ij} \; , \; 1\leq i,j\leq k,
\end{align*}
where $\delta_{ij}$ is the Kronecker delta. The above minimization problem is equivalent to the eigenvalue problem \cite{Kunisch01}
\begin{equation}\label{eg1}
\mathcal{W}\mathcal{W}^T\bm{M}\Psi_{i}=\sigma_{i}^2\Psi_{i} \; , \quad i=1,2,\ldots ,k,
\end{equation}
where $\mathcal{W}:= [\bm{w}_1 , \ldots, \bm{w}_{J}]\in\mathbb{R}^{N\times J}$ is the snapshot matrix whose $n$th column vector $\bm{w}_n$ is the solution vector of the FOM at time $t_n$. Then, the matrix $\bm{\Psi}$ of POD modes can be computed through the SVD of the snapshot matrix $\mathcal{W}$. In addition, between the coefficient vectors $\bm{w}(t)$ of the high-fidelity solutions and the coefficient vectors $\bm{w}^r(t)$ of reduced solutions, we have the relation
\begin{align*}
\bm{w}(t) \approx \bm{\Psi} \bm{w}^r(t), \quad \bm{w}^r(t) \approx \bm{\Psi}^T{\bm M} \bm{w}(t),
\end{align*}
from where we can find the initial reduced vector as $\bm{w}^r(0)=\bm{\Psi}^T{\bm M} \bm{w}_0$.  For a more detailed description, we refer to the study \cite{Karasozen17a}.

Projecting the system \eqref{semi_dg_rom} onto the reduced space $W_h^r$, and by the use of above formulations, we finally obtain the following ROM as a dynamical system:
\begin{equation}\label{rom}
  \frac{d}{dt} \bm{w}^r  =  -\bm{A}^{r} \bm{w}^r +  \bm{b}^{r}( \bm{w}^r),
\end{equation}
with the reduced stiffness matrix and the reduced nonlinear vector
\begin{align*}
\bm{A}^{r} &:= \bm{\Psi}^T\bm{A}\bm{\Psi}\in\mathbb{R}^{k\times k}\; , & (\bm{A}^{r})_{ij} = a_h(\psi_j(\bm{x}),\psi_i(\bm{x})),\\
\bm{b}^{r}( \bm{w}^r) &:= \bm{\Psi}^T\bm{b}( \bm{\Psi}\bm{w}^r)\in\mathbb{R}^{k}\; , & (\bm{b}^{r}( \bm{w}^r))_i = b_h(w_h^r(t,\bm{x}),\psi_i(\bm{x})),
\end{align*}
\begin{align*}
 b_h(w_h^r(t,\bm{x}),\psi_i(\bm{x})) = \sum \limits_{K \in \mathcal{T}_h} \int \limits_{K} f(w_h^r(t,\bm{x}))\psi_i(\bm{x}) \; d\bm{x}.
\end{align*}
Like the FOM \eqref{fom}, the ROM \eqref{rom} is solved in time by the AVF method.

\subsection{Approximation of the nonlinearities}

Although the reduced model \eqref{rom} is of small dimension, the computation of the nonlinear vector $\bm{b}^r( \bm{w}^r)$ still depends on the dimension $N$ of the FOM. In this section, we explain DEIM \cite{chaturantabut10nmr} to reduce the computational complexity due to the nonlinear vector in the ROM \eqref{rom}.

The DEIM aims  to find an approximation to the nonlinear vector $\bm{b}( \bm{\Psi}\bm{w}^r(t))$ which is the full dimensional nonlinear part of the reduced nonlinear vector $\bm{b}^r(\bm{w}^r)=\bm{\Psi}^T\bm{b}( \bm{\Psi}\bm{w}^r)$ for a.e. $t\in(0,T]$. It is followed by the  projection onto a subspace of the column space of the (nonlinear) snapshot matrix $\mathcal{B}:= [\bm{b}(\bm{w}_1) , \ldots, \bm{b}(\bm{w}_J)]\in\mathbb{R}^{N\times J}$ whose columns are the nonlinear vectors at discrete times $t_n$ for $n=1,\ldots , J$. Then, a DEIM basis $\{ \bm{Q}_i\}_{i=1}^m$ of dimension $m\ll N$ is constructed by applying the POD to the matrix $\mathcal{B}$, and the following approximation is used:
\begin{equation}\label{podG}
 \bm{b}( \bm{\Psi}\bm{w}^r(t)) \approx \bm{\hat{b}}( \bm{\Psi}\bm{w}^r(t)) = \bm{Q}\bm{c}(t),
\end{equation}
where $\bm{Q}:= [\bm{Q}_1\; \ldots \; \bm{Q}_m]\in\mathbb{R}^{N\times m}$ is the DEIM basis matrix, and $\bm{c}(t)$ is the unknown coefficients vector to be determined. Since the system \eqref{podG} is overdetermined, a projection matrix $\bm{P}=[e_{\mathfrak{p}_1},\ldots , e_{\mathfrak{p}_m}]\in\mathbb{R}^{N\times m}$ where $e_{\mathfrak{p}_i}$ is the $i$th column of the identity matrix $\bm{I}\in\mathbb{R}^{N\times N}$ is computed. Then, the approximate nonlinearity is derived as
\begin{equation}\label{deim_approx}
\bm{\hat{b}}( \bm{\Psi}\bm{w}^r(t)) =\bm{Q}(\bm{P}^T\bm{Q})^{-1}\bm{P}^T \bm{b}(\bm{\Psi} \bm{w}^r(t)),
\end{equation}
and then the reduced model \eqref{rom} can be rewritten as:
\begin{equation}\label{deim}
\frac{d}{dt} \bm{w}^r  =  -\bm{A}^r \bm{w}^r +  \bm{B}\bm{b}^r_{\text{deim}}( \bm{w}^r)
\end{equation}
where the matrix $\bm{B}:= \bm{\Psi}^T \bm{Q}(\bm{P}^T\bm{Q})^{-1}$ is computed once in the off-line stage, and the reduced nonlinear vector $\bm{b}^r_{\text{deim}}( \bm{w}^r):= \bm{P}^T \bm{b}(\bm{\Psi} \bm{w}^r)$ requires only $m\ll N$ integral evaluations. On the other hand, the computation of the Jacobian of the nonlinear vector requires $N\times n_p$ integral evaluations without DEIM, but it is only $m\times n_p$ with DEIM approximation.

For the details of the computation of the reduced non-linear vectors, we refer to the greedy DEIM algorithm  \cite{chaturantabut10nmr}.
Since the AVF method is an implicit time integrator, at each time step, a non-linear system of equations has to be solved by Newton's method. The reduced Jacobian has a block diagonal structure for the SIPG discretization, which is easily invertible \cite{Karasozen15}, and requires only $O(n_qN)$ operations with DEIM.

\subsection{ROM systems for the representative examples}

In this section, we give the ROM formulations related to the FOM systems \eqref{fom_rgle} and \eqref{fom_she} of RGL and SH equations, respectively. In all these two systems, the unknown solution is $w=(u,v)^T$ having two components. So, for a.e. $t\in (0,T]$, let us consider the dG (high-fidelity) solutions $u_h(t,\bm{x})\in W_h$ and $v_h(t,\bm{x})\in W_h$ with the related FOM solutions (coefficient vectors of the solutions) $\bm{u}(t)\in \mathbb{R}^N$ and $\bm{v}(t)\in \mathbb{R}^N$. Then, denote the corresponding reduced coefficient vectors by $\bm{u}^r(t)\in \mathbb{R}^{k_u}$ and $\bm{v}^r(t)\in \mathbb{R}^{k_v}$, which satisfy
\begin{align*}
\bm{u}(t) \approx \bm{\Psi_u} \bm{u}^r(t) \; &, \quad \bm{u}^r(t) \approx \bm{\Psi_u}^T{\bm M} \bm{u}(t), \\
\bm{v}(t) \approx \bm{\Psi_v} \bm{v}^r(t) \; &, \quad \bm{v}^r(t) \approx \bm{\Psi_v}^T{\bm M} \bm{v}(t), 
\end{align*}
where $\bm{\Psi_u}\in \mathbb{R}^{N\times{k_u}}$ and $\bm{\Psi_v}\in \mathbb{R}^{N\times{k_v}}$ are the matrices of POD modes related to the snapshot matrices
\begin{align*}
\mathcal{U} &:= [\bm{u}_1 , \ldots, \bm{u}_{J}]\in\mathbb{R}^{N\times J},\\
\mathcal{V} &:= [\bm{v}_1 , \ldots, \bm{v}_{J}]\in\mathbb{R}^{N\times J},
\end{align*}
of coefficient vectors of the high-fidelity solutions $u_h(t,\bm{x})$ and $v_h(t,\bm{x})$, respectively. Thus, the columns of the matrices of POD modes $\bm{\Psi_u}$ and $\bm{\Psi_v}$ are the coefficient vectors of the reduced basis functions $\psi_u(\bm{x})$ and $\psi_v(\bm{x})$ for the high-fidelity solutions $u_h(t,\bm{x})$ and $v_h(t,\bm{x})$, respectively, and they are the solutions of the eigenvalue problems \eqref{eg1} with $\mathcal{W}\in\{\mathcal{U},\mathcal{V}\}$, $\Psi_i\in\{\Psi_{u,i},\Psi_{v,i}\}$, and $k\in\{k_u,k_v\}$. Moreover, let $\bm{Q}_1\in\mathbb{R}^{N\times m_1}$ and $\bm{Q}_2\in\mathbb{R}^{N\times m_2}$ denote the DEIM basis of dimension $m_1$ and $m_2$, respectively, of the following (nonlinear) snapshot matrices:
\begin{align*}
\mathcal{B}_1 &:= [\bm{b}_1(\bm{u}_1,\bm{v}_1) ,\ldots, \bm{b}_1(\bm{u}_J,\bm{v}_J)]\in\mathbb{R}^{N\times J}, \\
\mathcal{B}_2 &:= [\bm{b}_2(\bm{u}_1,\bm{v}_1) ,\ldots, \bm{b}_2(\bm{u}_J,\bm{v}_J)]\in\mathbb{R}^{N\times J},
\end{align*}
with the corresponding DEIM projection matrices $\bm{P}_1$ and $\bm{P}_2$. Under these assumptions, the ROM with DEIM formulations for the RGL and SH equations are given as the following systems:\\

\noindent {\em ROM for RGL system \eqref{fom_rgle}\/}
\begin{equation}\label{rom_rgle}
\begin{aligned}
\frac{d}{dt}\bm{u}^r &= -\bm{A}^r_1\bm{u}^r + \bm{B}_1\bm{b}^r_{1,deim}(\bm{u}^r,\bm{v}^r), & \bm{b}^r_{1,deim}(\bm{u}^r,\bm{v}^r)= \bm{P}_1^T\bm{b}_1(\bm{u}^r,\bm{v}^r)\in\mathbb{R}^{m_1}, \\
\frac{d}{dt}\bm{v}^r &= -\bm{A}^r_2\bm{v}^r + \bm{B}_2\bm{b}^r_{2,deim}(\bm{u}^r,\bm{v}^r), & \bm{b}^r_{2,deim}(\bm{u}^r,\bm{v}^r)= \bm{P}_2^T\bm{b}_2(\bm{u}^r,\bm{v}^r)\in\mathbb{R}^{m_2},\\
& \bm{A}^r_1:= \bm{\Psi}_u^T\bm{A}\bm{\Psi}_u\in\mathbb{R}^{k_u\times k_u} , & \bm{B}_1:= \bm{\Psi}_u^T\bm{Q}_1(\bm{P}_1^T\bm{Q}_1)^{-1}\in\mathbb{R}^{k_u\times m_1},\\
& \bm{A}^r_2:= \bm{\Psi}_v^T\bm{A}\bm{\Psi}_v\in\mathbb{R}^{k_v\times k_v} , & \bm{B}_2:= \bm{\Psi}_v^T\bm{Q}_2(\bm{P}_2^T\bm{Q}_2)^{-1}\in\mathbb{R}^{k_v\times m_2}.
\end{aligned}
\end{equation}

\noindent {\em ROM for SH system \eqref{fom_she}\/}
\begin{equation}\label{rom_she}
\begin{aligned}
\frac{d}{dt}\bm{u}^r &= \bm{A}^r_1\bm{v}^r + \bm{B}_1\bm{b}^r_{1,deim}(\bm{u}^r,\bm{v}^r), & \bm{b}^r_{1,deim}(\bm{u}^r,\bm{v}^r)= \bm{P}_1^T\bm{b}_1(\bm{u}^r,\bm{v}^r)\in\mathbb{R}^{m_1}, \\
0 &= -\bm{A}^r_2\bm{u}^r + \bm{M}^r_1\bm{u}^r - \bm{M}^r_2\bm{v}^r, &  \bm{B}_1:= \bm{\Psi}_u^T\bm{Q}_1(\bm{P}_1^T\bm{Q}_1)^{-1}\in\mathbb{R}^{k_u\times m_1},\\
& \bm{A}^r_1:= \bm{\Psi}_u^T\bm{A}\bm{\Psi}_v\in\mathbb{R}^{k_u\times k_v} , & \bm{M}^r_1=\bm{\Psi}_u^T\bm{M}\bm{\Psi}_u\in\mathbb{R}^{k_u\times k_u},\\
& \bm{A}^r_2:= \bm{\Psi}_u^T\bm{A}\bm{\Psi}_u\in\mathbb{R}^{k_u\times k_u}, &  \bm{M}^r_2=\bm{\Psi}_u^T\bm{M}\bm{\Psi}_v\in\mathbb{R}^{k_u\times k_v}.
\end{aligned}
\end{equation}

\section{Energy analysis of ROM with DEIM approximation}
\label{Sec:analysis}

In this section, we investigate the gradient dissipative structures of reduced order models for RGL and SH systems based on DEIM approximation. We examine the difference of energy functionals at two successive time steps, i.e., $E_h[u^{r}_{n+1}, v^{r}_{n+1}]- E_h[u^{r}_{n}, v^{r}_{n}]$, and we find some upper bounds for it using the corresponding full discrete ROMs. The key point in the analysis is to justify the equivalence of the integrals of $f_1(u,v)$ and $f_2(u,v)$, which appear in the full discrete weak form due to AVF method, and the potential $F(u,v)$. Moreover, we include the contribution of DEIM by adding and subtracting the nonlinear terms to the difference of energy functionals written in terms of the weak forms, since DEIM offers an approximation for the nonlinear terms instead of exactly evaluating them. Then, we can examine the effect of DEIM on dissipative structures. On the other hand, the contribution of the parabolic-parabolic and parabolic-elliptic forms is reflected to the energy analysis. In detail, energy decreasing property of the gradient systems are proven under some condition on time step-size and this condition is related to both of the solutions of the parabolic-parabolic system, while the bound on time step-size depends only on the solution of the parabolic PDE for the parabolic-elliptic system. Moreover, we note that in the preceding analysis, the AVF method applied to linear terms coincides with the mid-point rule.

\subsection{Real Ginzburg-Landau system}

Related to the FOM \eqref{fom_rgle} of the RGL equation, we start from the reduced variational formulation, i.e. the variational formulation projected onto the reduced space $W_h^r$, and we apply the AVF method to obtain
\begin{subequations}\label{RGL_energy_1}
	\begin{align}
	\frac{1}{\Delta t}\langle u^{r}_{n+1}-u^r_{n}, \phi \rangle
	&=  - \frac{1}{2}a_h( u^{r}_{n+1}+u^{r}_{n},\phi)
	+ \frac{\mu}{2}\langle u^{r}_{n+1}+u^r_{n}, \phi \rangle + \int_0^1 \langle \tilde{f}_1(\tilde{u}^{r}, \tilde{v}^{r}),\phi\rangle d\xi,\\
	\frac{1}{\Delta t}\langle v^{r}_{n+1}-v^{r}_{n}, \eta \rangle
	&= - \frac{1}{2}a_h( v^{r}_{n+1}+v^{r}_{n},\eta)
	+ \frac{\mu}{2}\langle v^{r}_{n+1}+v^r_{n}, \eta \rangle + \int_0^1 \langle \tilde{f}_2(\tilde{u}^{r}, \tilde{v}^{r}),\eta\rangle d\xi,
	\end{align}
\end{subequations}
for all $\phi , \eta\in W_h^r$, where for analysis and computational purposes, we exclude the linear terms $\mu u$ and $\mu v$ from the nonlinearities, and we have set
\begin{align*}
& \tilde{f}_1(u,v) = - (u^2+v^2)u \; , & \tilde{f}_2(u,v) = - (u^2+v^2)v, \\
& \tilde{u}^{r} = \xi u_{n+1}^r + (1-\xi )u_{n}^r  \; , & \tilde{v}^{r} = \xi v_{n+1}^r + (1-\xi )v_{n}^r.
\end{align*}

Choosing the test functions as $\phi = u_{n+1}^r-u_{n}^r$ and $\eta = v_{n+1}^r - v_{n}^r$ from the reduced space $W_h^r$, and using the identity $(a+b)(a-b)=a^2-b^2$ on dG bilinear forms and linear terms on the right hand sides, we get
\begin{equation}\label{RGL_energy_2}
	\begin{aligned}
	\frac{1}{\Delta t}\| u_{n+1}^r-u_{n}^r\|^2 &= - \frac{1}{2}a_h( u_{n+1}^r,u_{n+1}^r) + \frac{1}{2}a_h( u_{n}^r,u_{n}^r) + \frac{\mu}{2}\| u^{r}_{n+1}\|^2 - \frac{\mu}{2}\| u^{r}_{n}\|^2 \\
	& + \int_0^1 \langle \tilde{f}_1(\tilde{u}^{r}, \tilde{v}^{r}), u_{n+1}^r-u_{n}^r \rangle d\xi, \\
	\frac{1}{\Delta t} \| v_{n+1}^r-v_{n}^r\|^2 &= - \frac{1}{2}a_h( v_{n+1}^r,v_{n+1}^r) + \frac{1}{2}a_h( v_{n}^r,v_{n}^r) + \frac{\mu}{2}\| v^{r}_{n+1}\|^2 - \frac{\mu}{2}\| v^{r}_{n}\|^2e \\
	& + \int_0^1 \langle \tilde{f}_2(\tilde{u}^{r}, \tilde{v}^{r}),v_{n+1}^r - v_{n}^r \rangle d\xi.
	\end{aligned}
\end{equation}
Now, we use the change of variables
\begin{align*}
z &= (\tilde{u}^r)^2+(\tilde{v}^r)^2, \\
dz &= 2( \tilde{u}^r(u_{n+1}^r-u_{n}^r) + \tilde{v}(v_{n+1}^r-v_{n}^r) ) d\xi,
\end{align*}
to obtain
\begin{equation}\label{RGL_energy_3}
\begin{aligned}
&\int_0^1 \langle \tilde{f}_1(\tilde{u}^{r}, \tilde{v}^{r}), u_{n+1}^r-u_{n}^r \rangle d\xi
+ \int_0^1 \langle \tilde{f}_2(\tilde{u}^{r}, \tilde{v}^{r}),v_{n+1}^r - v_{n}^r \rangle d\xi \\
&= - \int_0^1 \langle ((\tilde{u}^r)^2+(\tilde{v}^r)^2) ( \tilde{u}^r(u_{n+1}^r-u_{n}^r) + \tilde{v}^r(v_{n+1}^r-v_{n}^r) ) , 1\rangle  d\xi \\
&= - \frac{1}{2}\int_{(u_{n}^r)^2+(v_{n}^r)^2}^{(u_{n+1}^r)^2+(v_{n+1}^r)^2}\langle z,1\rangle dz \\
&= - \frac{1}{4} \langle\left( (u_{n+1}^r)^2+(v_{n+1}^r)^2  \right)^2 , 1 \rangle + \frac{1}{4} \langle\left( (u_{n}^r)^2+(v_{n}^r)^2  \right)^2, 1\rangle.
\end{aligned}
\end{equation}

Through the use of \eqref{RGL_energy_3}, the side-by-side summation of the equations in \eqref{RGL_energy_2} gives rise to
\begin{equation}\label{RGL_energy_4}
\begin{aligned}
0 &< \frac{1}{\Delta t}( \| u_{n+1}^r-u_{n}^r\|^2 + \| u_{n+1}^r-u_{n}^r\|^2 ) \\
&= -\frac{1}{2}a_h( u_{n+1}^r,u_{n+1}^r) - \frac{1}{2}a_h( v_{n+1}^r,v_{n+1}^r) + \langle F(u_{n+1}^r,v_{n+1}^r) , 1 \rangle \\
&\quad + \frac{1}{2}a_h( u_{n}^r,u_{n}^r) + \frac{1}{2}a_h( v_{n}^r,v_{n}^r) - \langle F(u_{n}^r,v_{n}^r) , 1 \rangle \\
&= \; - E_h[u^{r}_{n+1},v^{r}_{n+1}] + E_h[u^{r}_{n},v^{r}_{n}],
\end{aligned}
\end{equation}
where we used the fact that
\begin{equation*}
\langle F(u,v),1\rangle = \frac{\mu}{2}\| u\|^2 + \frac{\mu}{2}\| v\|^2 - \frac{1}{4}\langle (u^2+v^2)^2,1\rangle.
\end{equation*}
Hence, the identity \eqref{RGL_energy_4} implies that $E_h[u^{r}_{n},v^{r}_{n}] - E_h[u^{r}_{n+1},v^{r}_{n+1}]>0$, and that  without the use of DEIM approximation for the nonlinearities, the discrete energy of the ROM for RGL equation decreases as the time progresses.

\begin{remark}\label{remark_rom}
The identity \eqref{RGL_energy_4} implies that not only the ROM without DEIM approximation preserves the discrete energy decreasing property exactly, but also FOM \eqref{fom_rgle} of RGL equation preserves exactly the decreasing structure of the discrete energy. To see this, it is enough just to replace the reduced solution space $W_h^r$ with the full solution space $W_h$, and hence $u_h^r$ and $v_h^r$ are replaced with $u_h$ and $v_h$, respectively, in the above analysis.
\end{remark}

In the case of the ROM \eqref{rom_rgle} of RGL equation with the use of DEIM approximation for the nonlinearities, we should proceed  further. Before investigating the energy functional, we particularly consider the integral term in \eqref{RGL_energy_2} containing the nonlinearity $\tilde{f}_1(\tilde{u}^r_h, \tilde{v}^r_h)$ to obtain the following relation
\begin{equation}\label{RGL_energy_deim}
\begin{aligned}
& \int_{0}^{1}\left[ \int_{\Omega} \tilde{f}_1(\tilde{u}^r_h, \tilde{v}^r_h) (u^{r}_{n+1}-u^{r}_{n})\; dx \right] d \xi  \\
&= \int_{0}^{1}\left[ \sum_{i=1}^{k} (u_{i, n+1}^r(t) - u_{i, n}^r(t))  \int_{\Omega} \tilde{f}_1(\tilde{u}^r_h, \tilde{v}^r_h) \psi_{i}(\bm{x})  dx\right] d \xi \\
&= \int_{0}^{1} (\bm{\Psi}_u (\bm{u}^r_{n+1} - \bm{u}^r_{n}))^T \bm{\tilde{b}}_1(\bm{\Psi}_u\bm{\tilde{u}}^r(t),\bm{\Psi}_v\bm{\tilde{v}}(t)) d \xi.
\end{aligned}
\end{equation}
The relation for the integral term containing the nonlinearity $\tilde{f}_2(\tilde{u}^r_h, \tilde{v}^r_h)$ can be found in a similar manner. In order to catch the contribution of DEIM approximations, we need to replace the nonlinear vectors $\bm{\tilde{b}}_1$ and $\bm{\tilde{b}}_2$ with their DEIM approximations $\bm{\hat{b}}_1$ and $\bm{\hat{b}}_2$, respectively, given in \eqref{podG}. Now, we apply the replacements of DEIM approximations, we add and subtract the terms $\int_{0}^{1}\left[ \int_{\Omega} \tilde{f}_1(\tilde{u}^r_h, \tilde{v}^r_h) (u^{r}_{n+1}-u^{r}_{n})\; dx \right] d \xi$ and $\int_{0}^{1}\left[ \int_{\Omega} \tilde{f}_2(\tilde{u}^r_h, \tilde{v}^r_h) (v^{r}_{n+1}-v^{r}_{n})\; dx \right] d \xi$ to the right-hand side of \eqref{RGL_energy_2}, and we obtain the following relation \cite[Sec.3.3]{Karasozen17}

\begin{equation}\label{RGL_energy_5}
\begin{aligned}
&E_h[u^{r}_{n+1},v^{r}_{n+1}] - E_h[u^{r}_{n},v^{r}_{n}] \\
&= -\frac{1}{\Delta t} ( \| u^{r}_{n+1}-u^{r}_{n}\|^2 + \| v^{r}_{n+1}-v^{r}_{n}\|^2) \\
&+ \int_{0}^{1} (\bm{\Psi}_u (\bm{u}^r_{n+1} - \bm{u}^r_{n}))^T ( \bm{\tilde{b}}_1(\bm{\Psi}_u\bm{\tilde{u}}^r(t),\bm{\Psi}_v\bm{\tilde{v}}(t)) - \bm{\hat{b}}_1(\bm{\Psi}_u\bm{\tilde{u}}^r(t),\bm{\Psi}_v\bm{\tilde{v}}(t))) d \xi \\
&+ \int_{0}^{1} (\bm{\Psi}_v (\bm{v}^r_{n+1} - \bm{v}^r_{n}))^T ( \bm{\tilde{b}}_2(\bm{\Psi}_u\bm{\tilde{u}}^r(t),\bm{\Psi}_v\bm{\tilde{v}}(t)) - \bm{\hat{b}}_2(\bm{\Psi}_u\bm{\tilde{u}}^r(t),\bm{\Psi}_v\bm{\tilde{v}}(t))) d \xi,
\end{aligned}
\end{equation}
where the last two rows on the right hand side are additional because of the DEIM contribution. In particular, with the use of the integral mean theorem and the a priori estimates in \cite{chaturantabut10nmr, Chaturantabut12, Heinkenschloss14} for the DEIM approximation, we further obtain
\begin{equation}\label{RGL_energy_6}
\begin{aligned}
&\int_{0}^{1} (\bm{\Psi}_u (\bm{u}^r_{n+1} - \bm{u}^r_{n}))^T ( \bm{\tilde{b}}_1(\bm{\Psi}_u\bm{\tilde{u}}^r(t),
\bm{\Psi}_v\bm{\tilde{v}}(t)) - \bm{\hat{b}}_1(\bm{\Psi}_u\bm{\tilde{u}}^r(t),\bm{\Psi}_v\bm{\tilde{v}}(t))) d \xi  \\
&= \langle \bm{\Psi}_u (\bm{u}^r_{n+1} - \bm{u}^r_{n}), \bm{\tilde{b}}_1(\bm{\Psi}_u\bm{\hat{u}}^r(t),\bm{\Psi}_v\bm{\hat{v}}(t)) - \bm{\hat{b}}_1(\bm{\Psi}_u\bm{\hat{u}}^r(t),\bm{\Psi}_v\bm{\hat{v}}(t)) \rangle \\
&\leq \| \bm{\Psi}_u (\bm{u}^r_{n+1} - \bm{u}^r_{n}) \| \| \bm{\tilde{b}}_1(\bm{\Psi}_u\bm{\hat{u}}^r(t),\bm{\Psi}_v\bm{\hat{v}}(t)) - \bm{\hat{b}}_1(\bm{\Psi}_u\bm{\hat{u}}^r(t),\bm{\Psi}_v\bm{\hat{v}}(t)) \| \\
&\leq  \|  u^r_{n+1} - u^r_{n} \| \| \bm{R}^{-1} \|_2 \| (\bm{P}_1^T\bm{Q}_1)^{-1} \|_2 \| (\bm{I} - \bm{Q}_1 \bm{Q}_1^T) \bm{\tilde{b}}_1(\bm{\Psi}_u\bm{\hat{u}}^r(t),\bm{\Psi}_v\bm{\hat{v}}(t)) \|_2,
\end{aligned}
\end{equation}
for some $\bm{\hat{u}}_n^r$ between $\bm{u}_n^r$ and $\bm{u}^r_{n+1}$, where the matrix $\bm{R}$ is the Cholesky factor of the mass matrix $\bm{M}$, and $\|\cdot\|_2$ denotes the usual Euclidean inner product on $\mathbb{R}^{N}$. Applying the same idea for $\bm{\tilde{b}}_2(\bm{\Psi}_u\bm{\hat{u}}^r(t),\bm{\Psi}_v\bm{\hat{v}}(t))$, combining \eqref{RGL_energy_5} and \eqref{RGL_energy_6}, we finally obtain
\begin{align}\label{RGL_energy_7}
&E_h[u^{r}_{n+1},v^{r}_{n+1}] - E_h[u^{r}_{n},v^{r}_{n}] \nonumber \\
&\leq \|u^{r}_{n+1}-u^{r}_{n}\|^2 \left(-\frac{1}{\Delta t}
+ \frac{\| \bm{R}^{-1}\| \|(\bm{P}_1^T\bm{Q}_1 )^{-1}|| \|(\bm{I} -\bm{Q}_1\bm{Q}^T_1) \bm{\tilde{b}}_1(\bm{\Psi}_u\bm{\hat{u}}^r(t),\bm{\Psi}_v\bm{\hat{v}}(t)) ||
	}{\|u^{r}_{n+1}-u^{r}_{n}\|} \right) \nonumber\\
&+ \| v^{r}_{n+1}-v^{r}_{n}\|^2 \left(-\frac{1}{\Delta t}
+ \frac{\| \bm{R}^{-1}\| \|(\bm{P}_2^T\bm{Q}_2 )^{-1}|| \|(\bm{I} -\bm{Q}_2\bm{Q}^T_2) \bm{\tilde{b}}_2(\bm{\Psi}_u\bm{\hat{u}}^r(t),\bm{\Psi}_v\bm{\hat{v}}(t)) ||
}{\|v^{r}_{n+1}-v^{r}_{n}\|} \right).
\end{align}

We deduce that ROM with DEIM formulation preserves the energy decreasing property, i.e., $E_h[u^{r}_{n+1},v^{r}_{n+1}] \leq  E_h[u^{r}_{n},v^{r}_{n}]$, if the time step-size satisfies
\begin{equation}\label{upperbound}
\begin{aligned}
\Delta t
\leq \frac{1}{\| \bm{R}^{-1}\| } \times
& \min \Big\lbrace \frac{ \|u^{r}_{n+1}-u^{r}_{n}\|}{\|(\bm{P}_1^T\bm{Q}_1 )^{-1}|| \|(\bm{I} -\bm{Q}_1\bm{Q}^T_1) \bm{\tilde{b}}_1(\bm{\Psi}_u\bm{\hat{u}}^r(t),\bm{\Psi}_v\bm{\hat{v}}(t)) || }, \\
& \qquad \frac{ \|v^{r}_{n+1}-v^{r}_{n}\|}{\|(\bm{P}_2^T\bm{Q}_2 )^{-1}|| \|(\bm{I} -\bm{Q}_2\bm{Q}^T_2) \bm{\tilde{b}}_2(\bm{\Psi}_u\bm{\hat{u}}^r(t),\bm{\Psi}_v\bm{\hat{v}}(t)) ||} \Big\rbrace.
\end{aligned}
\end{equation}

By construction, the matrices $\bm{Q}_1$ and $\bm{Q}_2$ are orthonormal, and the terms $\|(\bm{P}_1^T\bm{Q}_1 )^{-1}\|_2$ and $\|(\bm{P}_2^T\bm{Q}_2 )^{-1}\|_2$ are of moderate magnitude. They vary in the numerical tests between $10-30$, and $30-60$, respectively.  Hence, the upper bound for the time step-size $\Delta t$ in the right hand side of \eqref{upperbound} is a sufficiently large number so that the discrete energy for the ROM \eqref{rom_rgle} of RGL equation with the use of DEIM approximation decreases almost unconditionally for large time step-sizes.

\subsection{Swift-Hohenberg system}

The AVF applied reduced variational formulation related to the FOM \eqref{fom_she} of the SH equation is given by
\begin{equation}\label{SH_energy_1}
\begin{aligned}
\frac{1}{\Delta t} \langle u^{r}_{n+1} - u^{r}_{n}, \phi \rangle 
&= \frac{1}{2}a_h( v^{r}_{n+1}+v^{r}_{n},\phi) - \frac{1}{2}\langle v^{r}_{n+1}+v^{r}_{n} , \phi \rangle + \int_0^1 \langle f(\tilde{u}^{r}),\phi\rangle d\xi, \\
0 &= -\frac{1}{2}a_h(u^{r}_{n+1}+u^{r}_{n}, \eta) +  \frac{1}{2}\langle u^{r}_{n+1}+u^{r}_{n} , \eta \rangle - \frac{1}{2}\langle v^{r}_{n+1}+v^{r}_{n} , \eta \rangle,
\end{aligned}
\end{equation}
where, again for computational purposes as we did for RGL equation, we have set the nonlinearity $f(u)=\mu u - u^3$, and $\tilde{u}^{r} = \xi u^{r}_{n+1} + (1-\xi)u^{r}_{n}$. Choosing $\phi=u^{r}_{n+1}- u^{r}_{n}$ and $\eta=-(v^{r}_{n+1}- v^{r}_{n})$ in \eqref{SH_energy_1}, summing up the equations side-by-side, using the bi-linearity of $a_h(\cdot,\cdot)$, and from the identity $(a+b)(a-b)=a^2-b^2$, we get
\begin{equation}\label{SH_energy_2}
\begin{aligned}
\frac{1}{\Delta t} \|u^{r}_{n+1}-u^{r}_{n}\|^2 =& \; a_h( u^{r}_{n+1},v^{r}_{n+1}) - a_h( u^{r}_{n},v^{r}_{n}) + \frac{1}{2}\|v^{r}_{n+1}\|^2 - \frac{1}{2}\|v^{r}_{n}\|^2 \\
& - \langle u^{r}_{n+1} , v^{r}_{n+1} \rangle + \langle u^{r}_{n} , v^{r}_{n} \rangle + \int_0^1 \langle f(\tilde{u}^{r}),u^{r}_{n+1}-u^{r}_{n}\rangle d\xi
\end{aligned}
\end{equation}
With the change of variables
\begin{equation*}
z:= \tilde{u}^{r} = \xi u^{r}_{n+1} + (1-\xi)u^{r}_{n} \; , \qquad dz = (u^{r}_{n+1}-u^{r}_{n})d\xi,
\end{equation*}
we obtain that
\begin{equation}\label{SH_energy_3}
\begin{aligned}
& \int_0^1 \langle f(\tilde{u}^{r}),u^{r}_{n+1}-u^{r}_{n}\rangle d\xi = \int_0^1 \langle \mu \tilde{u}^{r} - (\tilde{u}^{r})^3, u^{r}_{n+1}-u^{r}_{n}\rangle d\xi \\
&= \int_{u_{n}^r}^{u_{n+1}^r}\langle \mu z - z^3,1\rangle dz \\
&= \frac{\mu}{2} \| u_{n+1}^r\|^2 - \frac{1}{4} \langle  (u_{n+1}^r)^4, 1\rangle - \frac{\mu}{2} \| u_{n}^r\|^2 + \frac{1}{4} \langle  (u_{n}^r)^4, 1\rangle. 
\end{aligned}
\end{equation}
Combining \eqref{SH_energy_2} and \eqref{SH_energy_3}, we obtain that
\begin{align*}
0 <& \frac{1}{\Delta t} \|u^{r}_{n+1}-u^{r}_{n}\|^2 \\
&= a_h(u_{n+1}^r,v_{n+1}^r) + \langle F(u_{n+1}^r,v_{n+1}^r),1\rangle - a_h(u_{n}^r,v_{n}^r) - \langle F(u_{n}^r,v_{n}^r),1\rangle \\
&= -E_h[u^{r}_{n+1}, v^{r}_{n+1}] + E_h[u^{r}_{n}, v^{r}_{n}],
\end{align*}
where we used the fact that
\begin{equation*}
\langle F(u,v),1\rangle = \frac{\mu}{2} \| u\|^2 - \frac{1}{4} \langle  u^4, 1\rangle - \langle u, v\rangle + \frac{1}{2} \| v\|^2.
\end{equation*}
Clearly, the ROM without DEIM approximation, and also the FOM for SH equation (Remark~\ref{remark_rom}), preserves the energy decreasing property exactly.

For the contribution of the DEIM approximation, we follow the similar steps as we did for the RGL equation: we replace the nonlinear vector with its DEIM approximation, we add and subtract the term $\int_{0}^{1}\langle f(\tilde{u}^r_h), u^{r}_{n+1}-u^{r}_{n}\rangle d \xi$ to the right-hand side of \eqref{SH_energy_3} to reach at
\begin{equation}\label{SH_energy_4}
\begin{aligned}
&E_h[u^{r}_{n+1}, v^{r}_{n+1}]- E_h[u^{r}_{n}, v^{r}_{n}] \\
&= \; -\frac{1}{\Delta t} \| u^{r}_{n+1}- u^{r}_{n} \|^2  + \langle \bm{\Psi}_u(\bm{u}^r_{n+1}-\bm{u}^r_{n}), \bm{b}(\bm{\Psi}_u\bm{\tilde{u}}^r) - \bm{\hat{b}}(\bm{\Psi}_u\bm{\tilde{u}}^r) \rangle, \\
&\leq\; -\frac{1}{\Delta t} \| u^{r}_{n+1}- u^{r}_{n} \|^2  \\
&+ \|  u^r_{n+1} - u^r_{n} \| \| \bm{R}^{-1} \|_2 \| (\bm{P}^T\bm{Q})^{-1} \|_2 \| (\bm{I} - \bm{Q} \bm{Q}^T) \bm{b}(\bm{\Psi}_u\bm{\hat{u}}^r) \|_2,
\end{aligned}
\end{equation}
for some $\bm{\hat{u}}_r^n$ between $\bm{u}_r^n$ and $\bm{u}_r^{n+1}$, and with the DEIM approximation $\bm{b}( \bm{\Psi}\bm{w}^r(t)) \approx  \bm{\hat{b}}( \bm{\Psi}\bm{w}^r(t))$. From \eqref{SH_energy_4}, we find out that the discrete energy decreasing property, i.e., $E_h[u^{r}_{n+1}, v^{r}_{n+1}] \leq E_h[u^{r}_{n}, v^{r}_{n}]$, holds for ROM \eqref{rom_she} of SH equation with the DEIM approximation under the following condition on the step-size
\begin{align}\label{SH_energy_5}
\Delta t
&\leq \frac{ \|u^{r}_{n+1}- u^{r}_{n}\|}{\| \bm{R}^{-1}\|_2 \|(\bm{P}^T\bm{Q} )^{-1}||_2 \|(\bm{I} -\bm{Q}\bm{Q}^T) \bm{b}(\bm{\Psi}_u\bm{\hat{u}}^r) \|_2},
\end{align}
which allows the energy decreasing property almost unconditionally, by the same reasoning we have stated for the RGL equation.

\section{Numerical results}
\label{Sec:numeric}

In this section, we present some numerical results for the two gradient systems. In all simulations, we use linear dG basis functions on a uniform $32\times 32$ rectangular grid with $2048$ triangular elements in the rectangular domain $\Omega =[a,b]^2$. The number of POD and DEIM modes are determined by a user-defined threshold $\epsilon$ such that "{\em relative information content\/}" criterion
\begin{equation}\label{ric}
RIC(k) = \frac{\sum_{i=1}^k \sigma_i^2}{\sum_{i=1}^{d_z} \sigma_i^2},
\end{equation}
satisfies that $RIC(k) \ge 1 - \varepsilon$. In \eqref{ric} $\sigma_i$ is the $i$th singular value, and $d_z$ is the rank of the aforementioned snapshot matrix.  When $\varepsilon = 0$, the snapshot data is recovered perfectly up to numerical errors. In the  numerical simulations, we took for the DEIM modes a smaller value of threshold $\epsilon$ than for the POD  so that the POD and DEIM errors have the same level of accuracy.  For the computation of the SVD, we use the randomized SVD (rSVD) which is computationally much faster then the standard SVD. 
For a matrix $A\in\mathbb{R}^{m\times n}$ together with a target rank $k$, the rSVD uses a random $n\times k$ matrix.
The approximation error is controlled by an oversampling  parameter $p$, and taking a random $n\times (k+p)$ matrix instead. More accurate approximation is obtained by increasing the number of oversampling parameter.  For the choice of the oversampling parameter  $p= k$, the following error bound is obtained for the rSVD \cite{Halko11}
    $$
    \mathbb{E} [\; ||A- U_l\Sigma_lV_l^*||\;] = \sigma_{l+1} \left [ 1 + 4\sqrt{\frac{2\min (m,n)}{l-1}} \right ]^{\frac{1}{2q +1}},
    $$
    where $q$ is the number of the power iterations in the rSVD, $l=2k$ and $\sigma_{l+1}$ is the $(l+1)$th singular value. We have used only one power iteration $q=1$.
Significant speedups over deterministic SVD algorithms are obtained for large matrices with rapidly decaying singular values.
Computational complexity for deterministic SVD algorithms is ${\mathcal O} (mn^2)$ when $m > n$. Computational cost of rSVD is  ${\mathcal O} (mnk)$ \cite{Erichson16}. For the details of rSVD algorithm we refer to
\cite{Alla16a, Karasozen18}.

In the numerical experiments, we have used the same time-step size $\Delta t =0.01$ for the FOMs and ROMs. Only one Newton iteration is needed in solving nonlinear equations of  the FOMs and ROMs.  All simulations are performed on a Windows 10 machine with Intel Core i7, 2.5 GHz and 8 GB using MATLAB R2014.  CPU time is measured in seconds.

The numerical simulations are stopped until  the system reaches a spatially inhomogeneous steady-state satisfying, for a user defined tolerance $tol$, the following criteria
$$
\frac{||u^{n+1} - u^n ||}{||u^n||} \le tol , \qquad \frac{||v^{n+1} - v^n ||}{||v^n||} \le tol.
$$

The errors between the  FOM and ROM solutions are given in $L^2$-$L^2$ norm (Frobenious norm) over the discrete space and time, whereas the energy errors are computed
in $L^2$ norm over the time interval. We present the FOM and ROM solutions for the GL equation only for the $u$ component, because the patterns for $v$ are similar to $u$.

\subsection{Ginzburg-Landau equation}

 The numerical example is taken from  \cite{Patra13} with $\mu = 0.5$ on a rectangular domain $\Omega = [0,256] \times [0,256]$ under periodic boundary conditions.
Initial conditions are chosen as random perturbations around a constant solution  $u({\mathbf x},0)=v({\mathbf x},0)=2 + {\mathrm {rand}}({\mathbf x})$, where ${\mathrm {rand}}: \Omega\rightarrow [0,1]$  is a uniform random distribution.  The steady state solutions are computed with $tol=10^{-4}$.

\begin{figure}[htb!]
\centering
\subfloat{\includegraphics[scale=0.6]{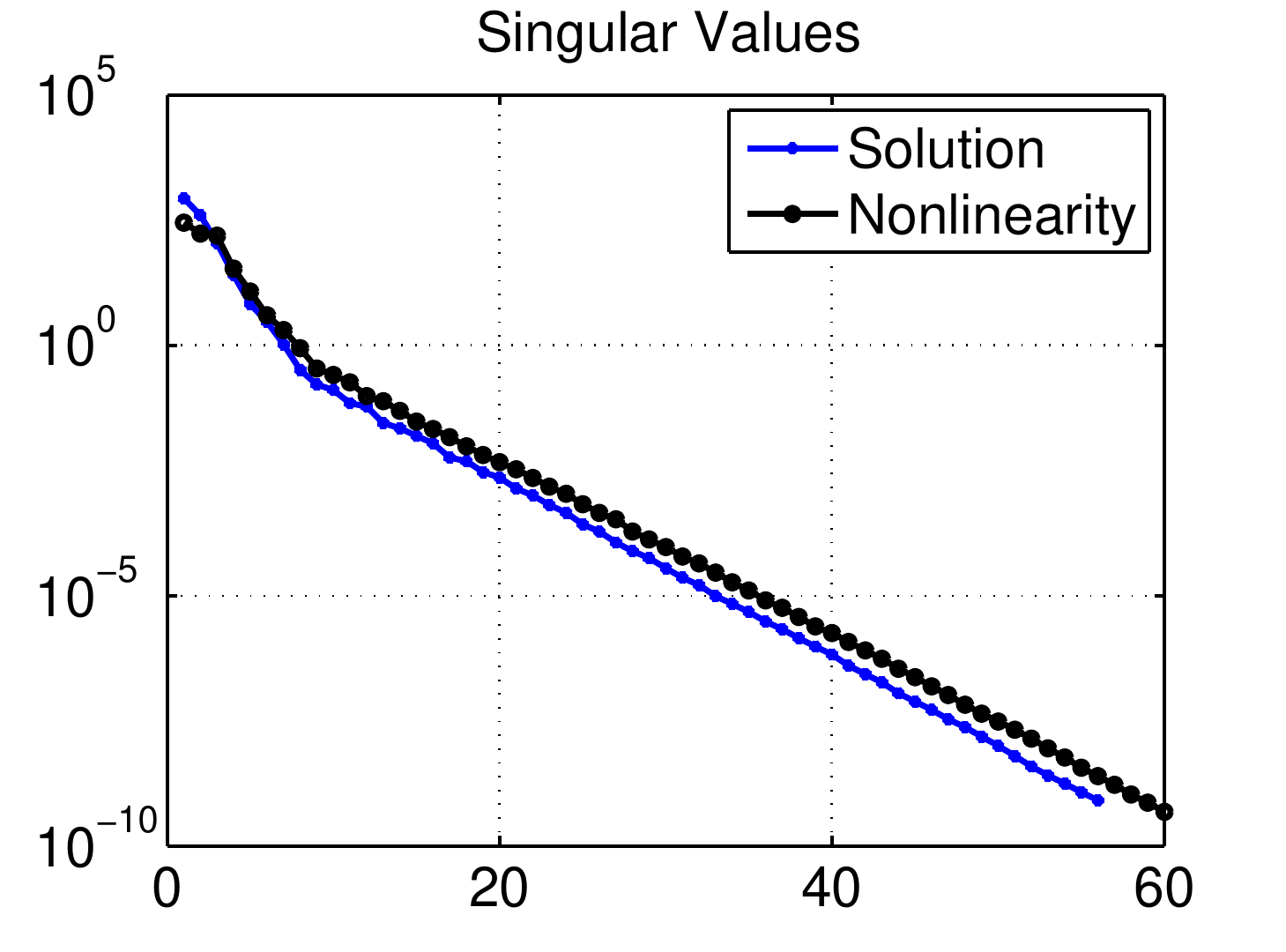}}
\caption{RGL: Decay of singular values for the GLE\label{gle_svd}}
\end{figure}

\begin{figure}[htb!]
\centering
\subfloat{\includegraphics[scale=0.33]{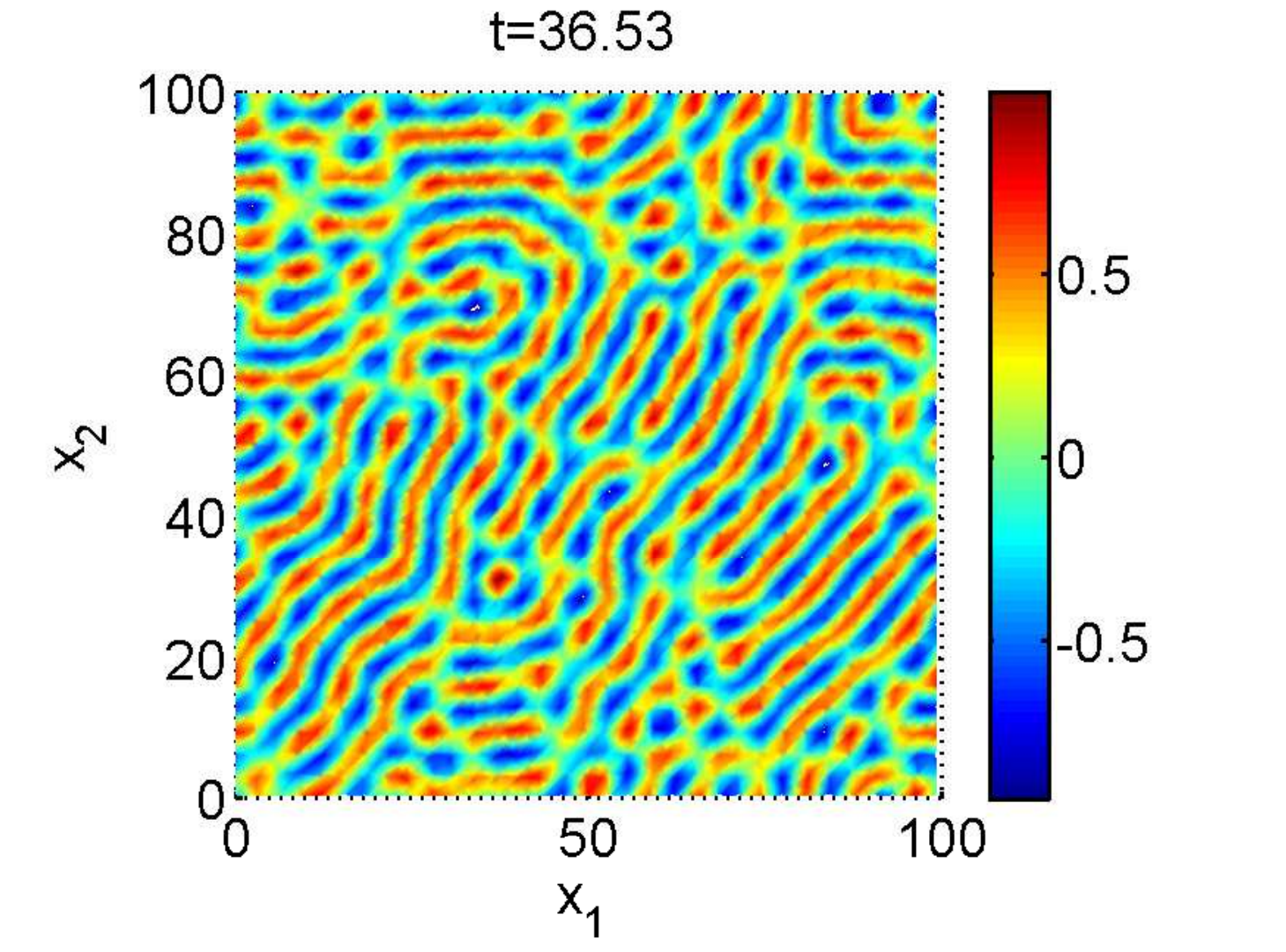}}
\subfloat{\includegraphics[scale=0.33]{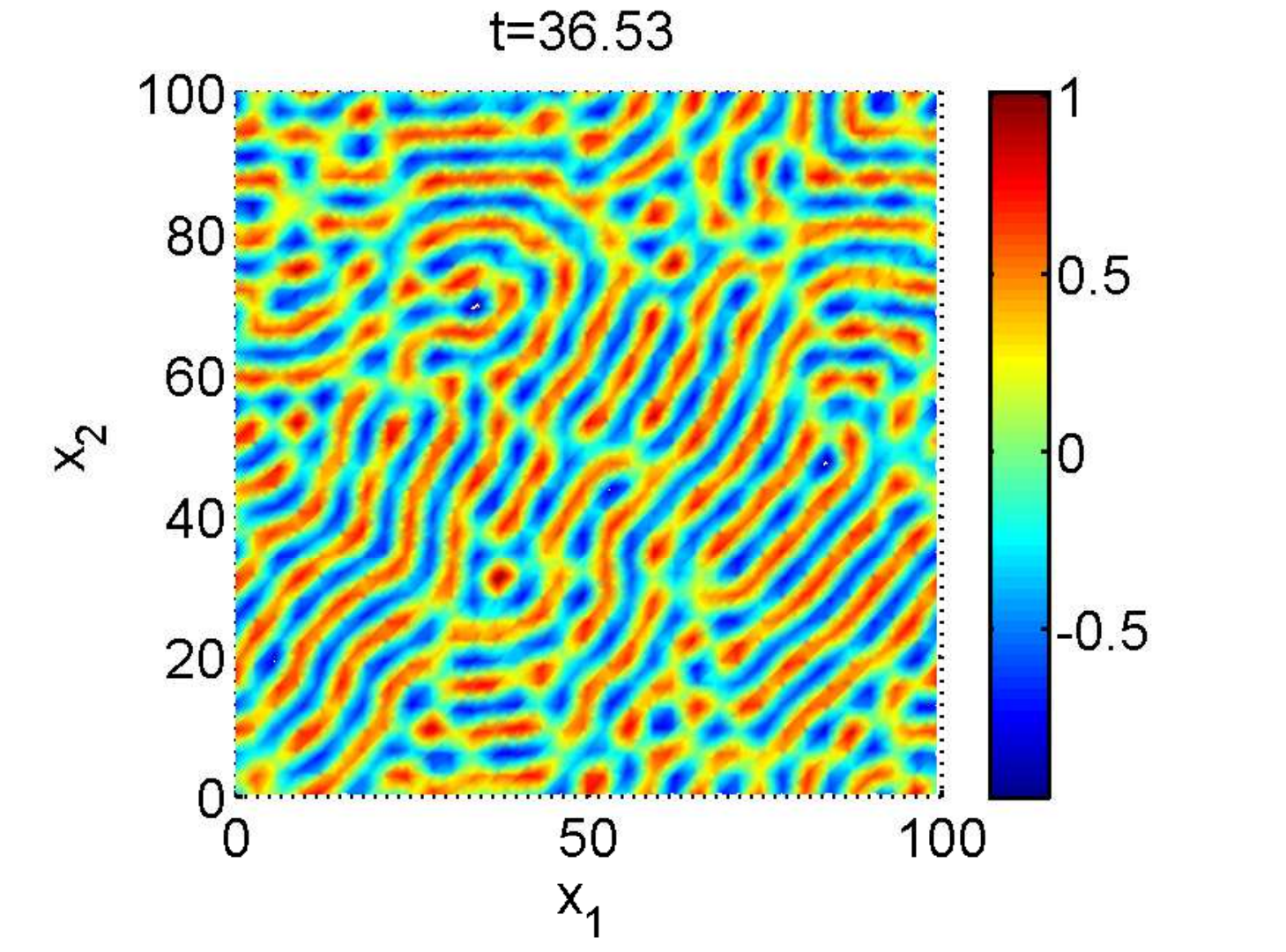}}
\subfloat{\includegraphics[scale=0.33]{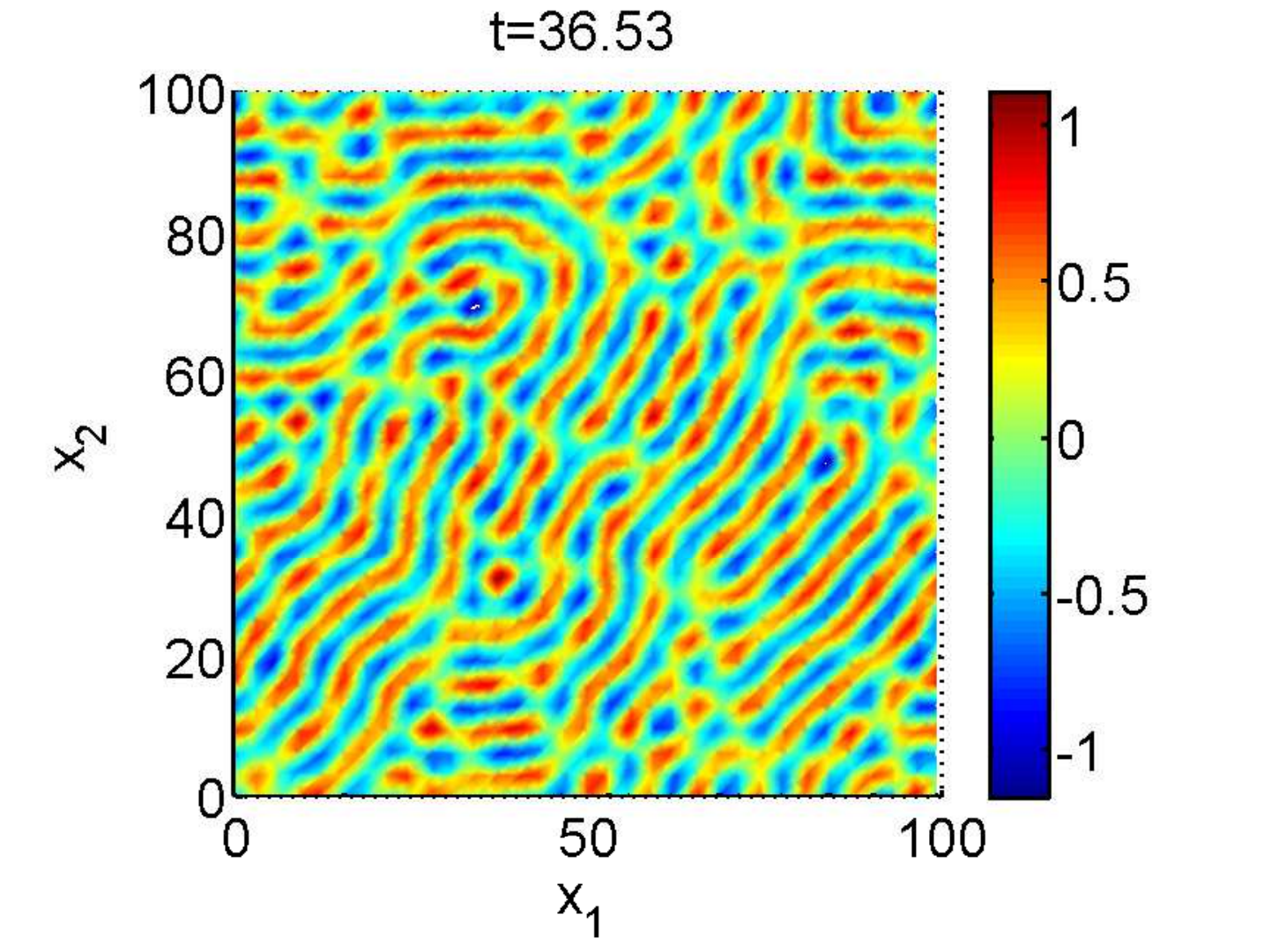}}
\caption{RGL: Steady state solutions: FOM, POD, DEIM\label{gle_rom}}
\end{figure}

\begin{figure}[htb!]
\centering
\subfloat{\includegraphics[scale=0.5]{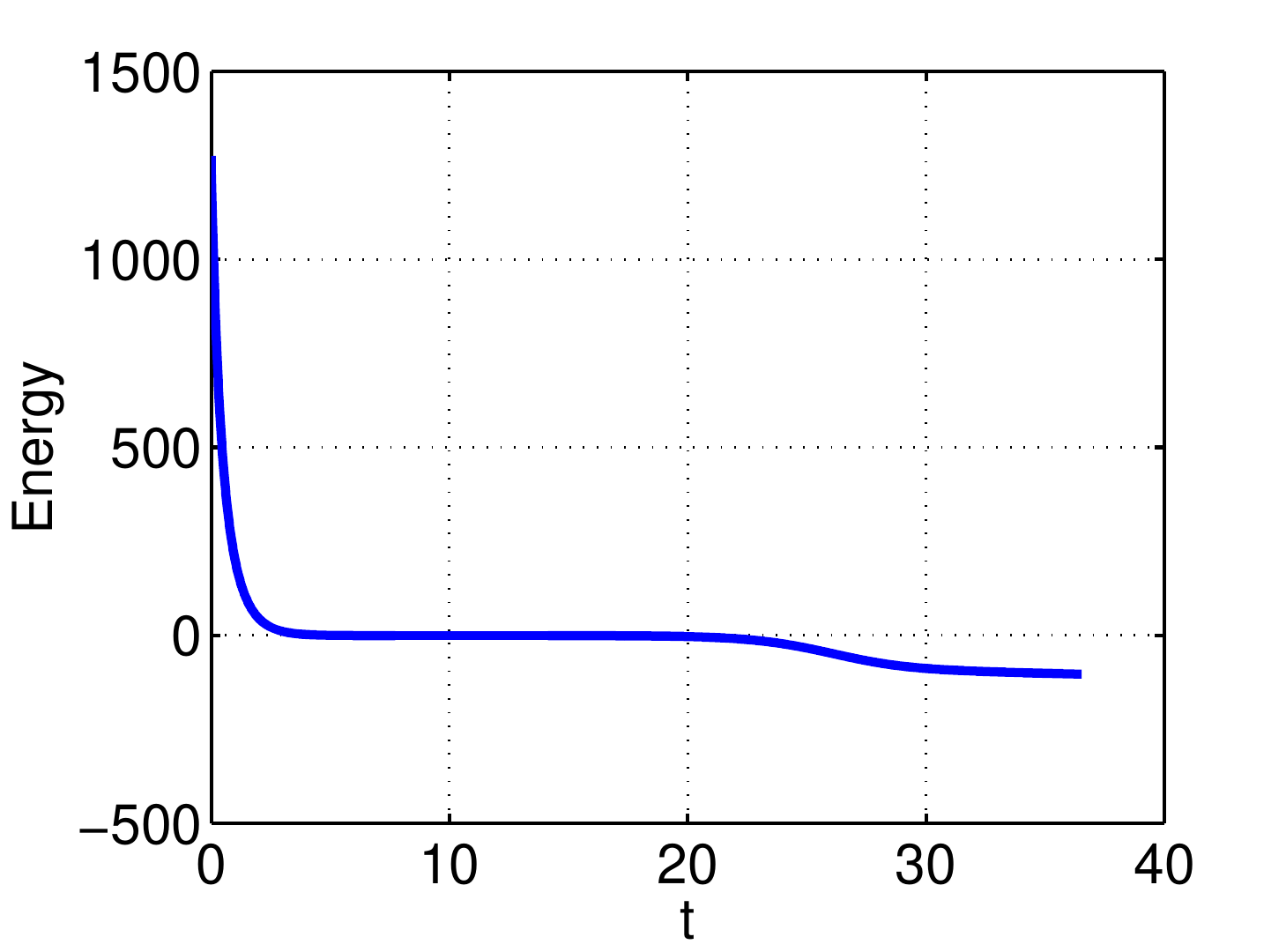}}
\subfloat{\includegraphics[scale=0.5]{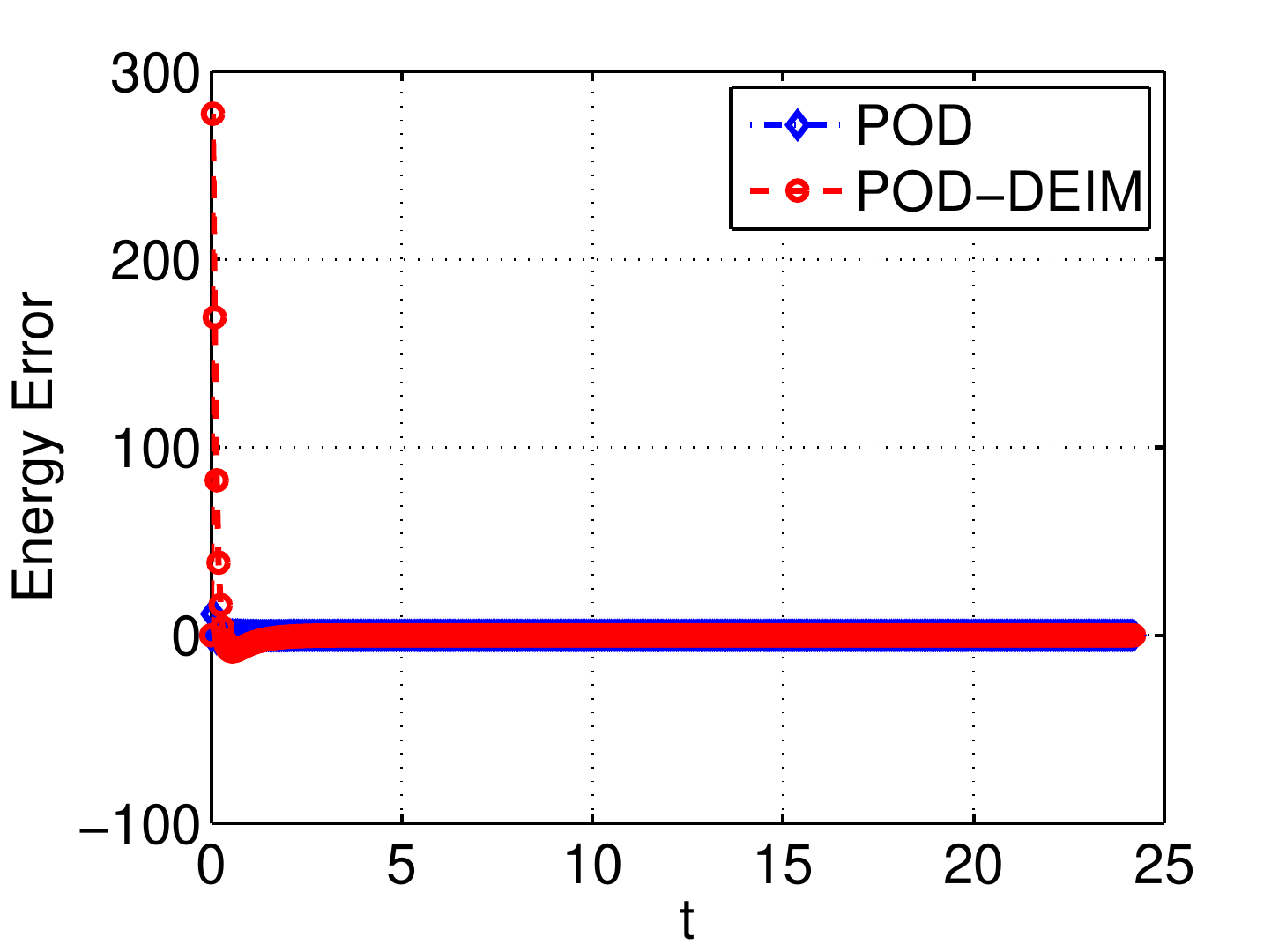}}
\caption{RGL: Energy decrease and energy error\label{gle_energy}}
\end{figure}

The singular values are plotted in semi-logarithmic form  in Fig. \ref{gle_svd}, with the index of singular values on the horizontal axis. The singular values of the snapshot matrices and of the nonlinear terms decay linearly in Fig \ref{gle_svd} . The threshold for the selection of POD and DEIM modes are chosen as $\epsilon = 10^{-4}$ and $\epsilon =10^{-6}$, respectively. In Fig. \ref{gle_rom} the FOM solutions for 4 POD and 5 DEIM modes are shown at the steady state. For a small number of POD and DEIM modes, the FOM and ROM are very close in Fig.~\ref{gle_rom} and the $L^2$-$L^2$ errors, given in Table~\ref{gle_error}, are relatively small. The FOM energy decays rapidly for a short time, and around $T=22$  it remains constant  until the steady-state in Fig.~\ref{gle_energy}. The energy is very well preserved by the ROMs over the whole time integration as shown in Fig.~\ref{gle_energy} and Table~\ref{gle_error}. Using the DEIM modes, the computational cost is reduced about a factor of 7 over the  POD solutions, see Table~\ref{gle_cpu}.

\begin{table}[H]
\centering
\caption{RGL: POD \& DEIM, energy  errors\label{gle_error}}
\begin{tabular}{llllll}
\hline
    \# POD modes  &     \#DEIM modes &      POD   &    DEIM  & POD & DEIM \\
 u(v)   &  u(v) & error &  error & energy error & energy error \\  \hline
       4(4)  &     5(5) &         5.75e-04 &   1.56e-03 & 1.81e-05 & 1.23e-04\\
\hline
\end{tabular}
\end{table}

\begin{table}[H]
\centering
\caption{RGL: CPU time and speedup factors\label{gle_cpu}}
\begin{tabular}{ l c c }
\hline
   & CPU time & speedup \\
\hline
FOM            & 1519.96 & -  \\
POD &  787.44 & 1.93 \\
DEIM  & 104.0   &   14.60 \\
\hline
\end{tabular}
\end{table}

\subsection{Swift-Hohenberg equation}
For numerical simulations, we choose the Example 5.5 in  \cite{Deghan17} with the parameter value $\mu = 0.3$  in the spatial  domain $\Omega = [0,100]\times[0,100]$ under periodic boundary conditions.
Initial condition is taken as random perturbation around the equilibrium  \cite{Deghan17,Moreno14} bounded below and above by $\mp 5\times 10^{-5}$. FOM solution is computed until the steady state with $tol =10^{-7}$. POD and DEIM modes are chosen with the threshold values $\epsilon = 10^{-6}$ and $\epsilon =10^{-8}$, respectively.

\begin{figure}[htb!]
\centering
\subfloat{\includegraphics[scale=0.6]{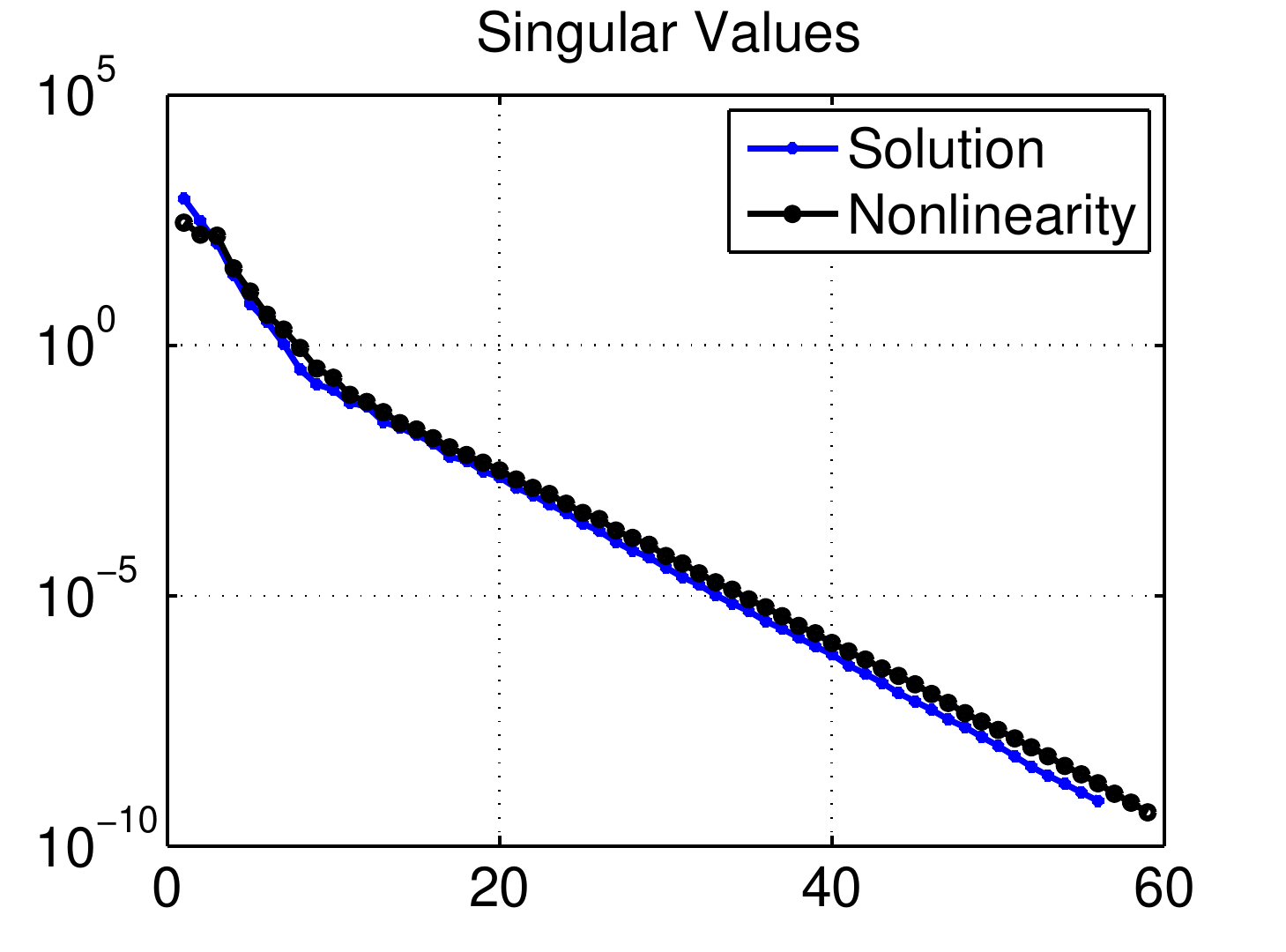}}
\caption{SH: Decay of singular values of SH equation\label{sh_svd}}
\end{figure}

\begin{figure}[htb!]
\centering
\subfloat{\includegraphics[scale=0.33]{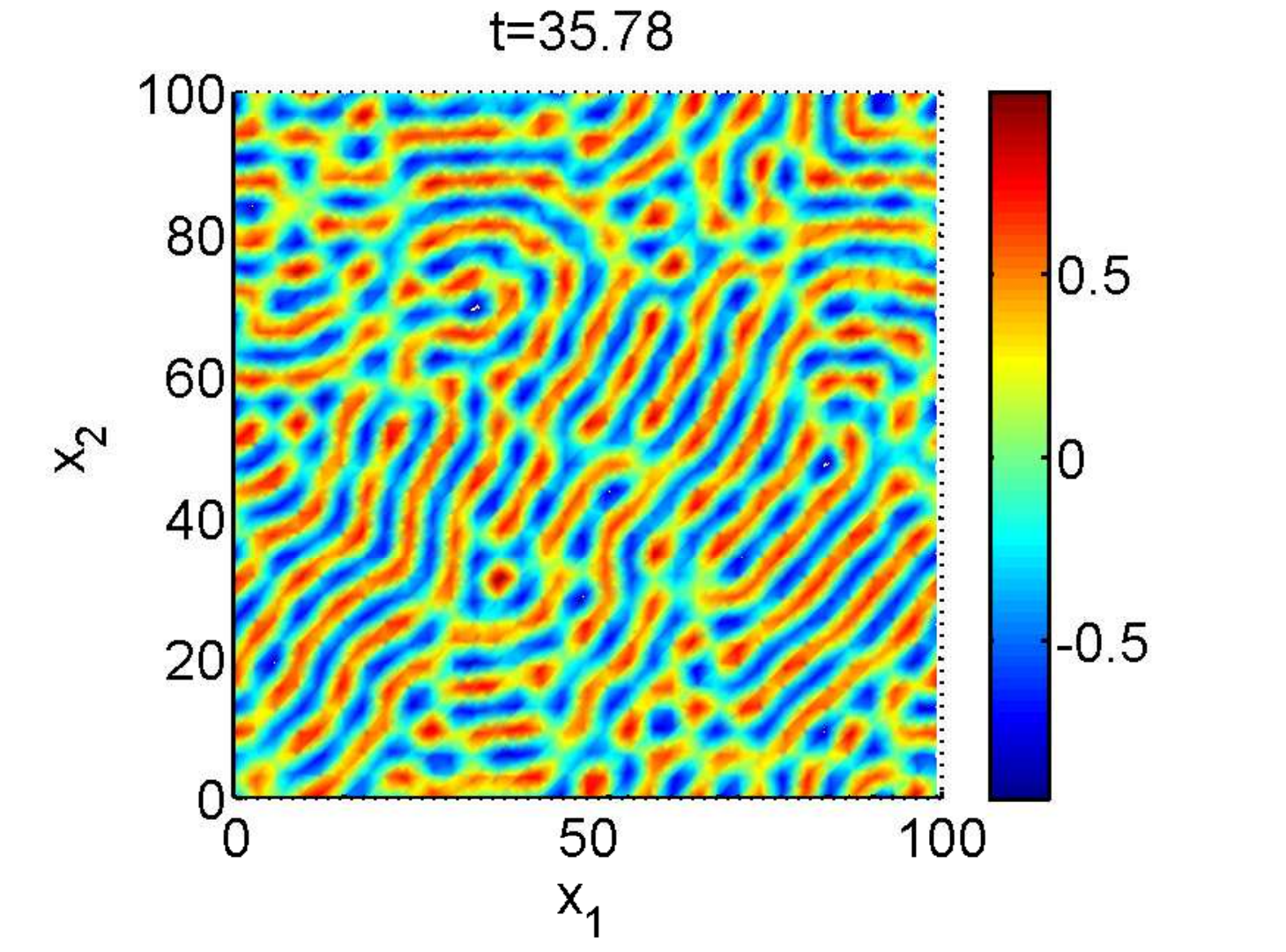}}
\subfloat{\includegraphics[scale=0.33]{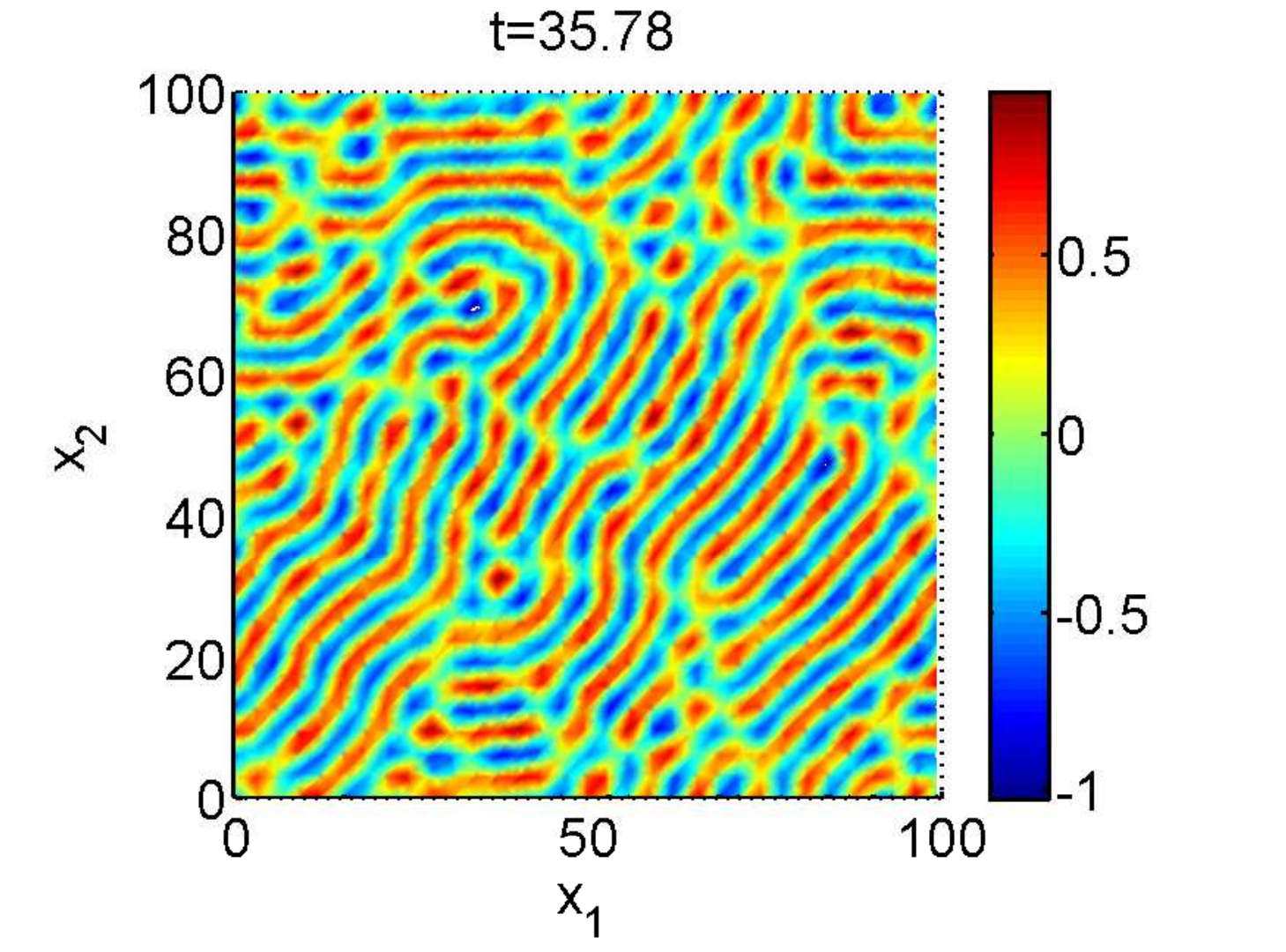}}
\subfloat{\includegraphics[scale=0.33]{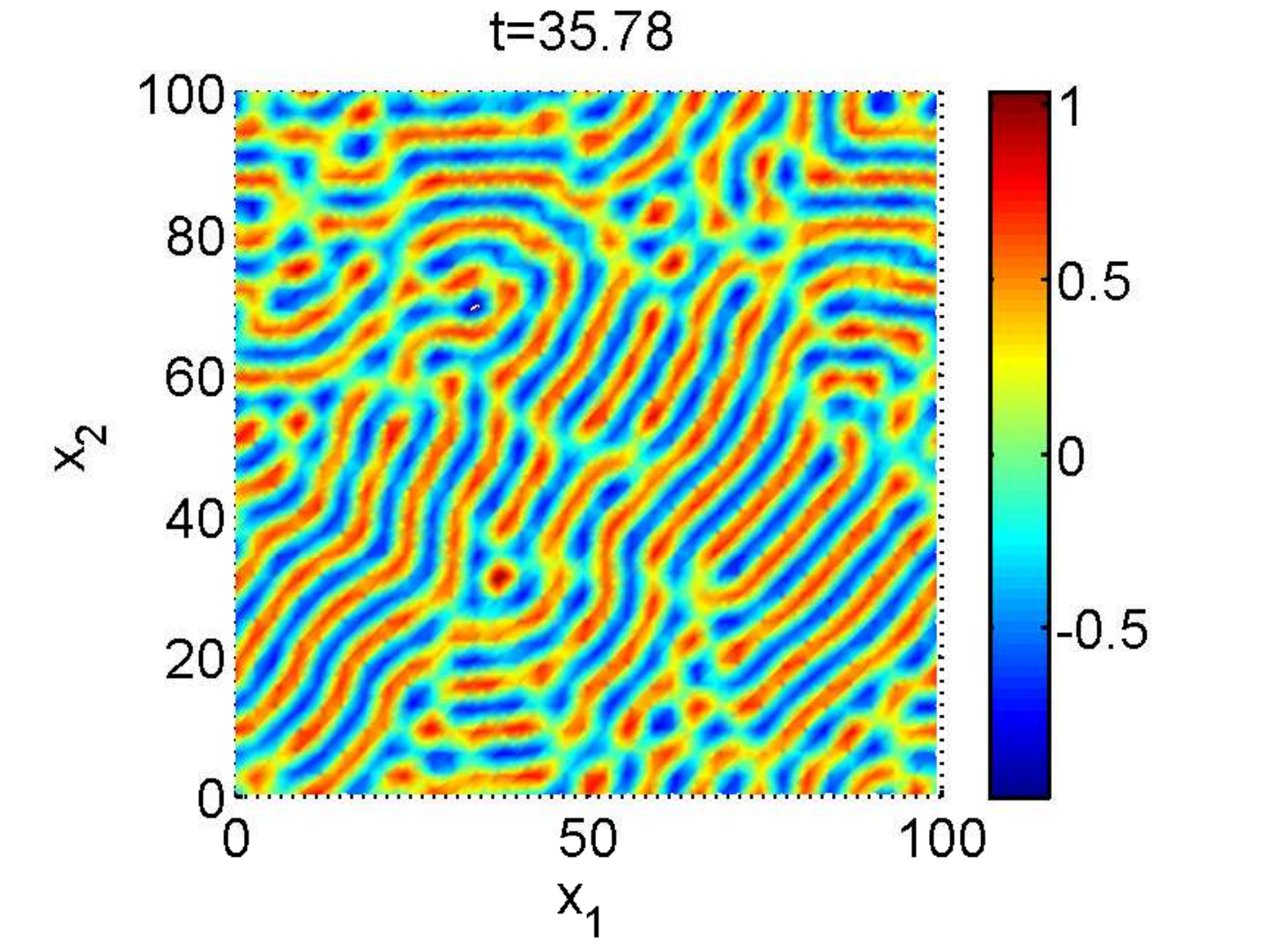}}
\caption{SH: Steady state solutions: FOM, POD, DEIM\label{sh_rom}}
\end{figure}

\begin{figure}[htb!]
\centering
\subfloat{\includegraphics[scale=0.5]{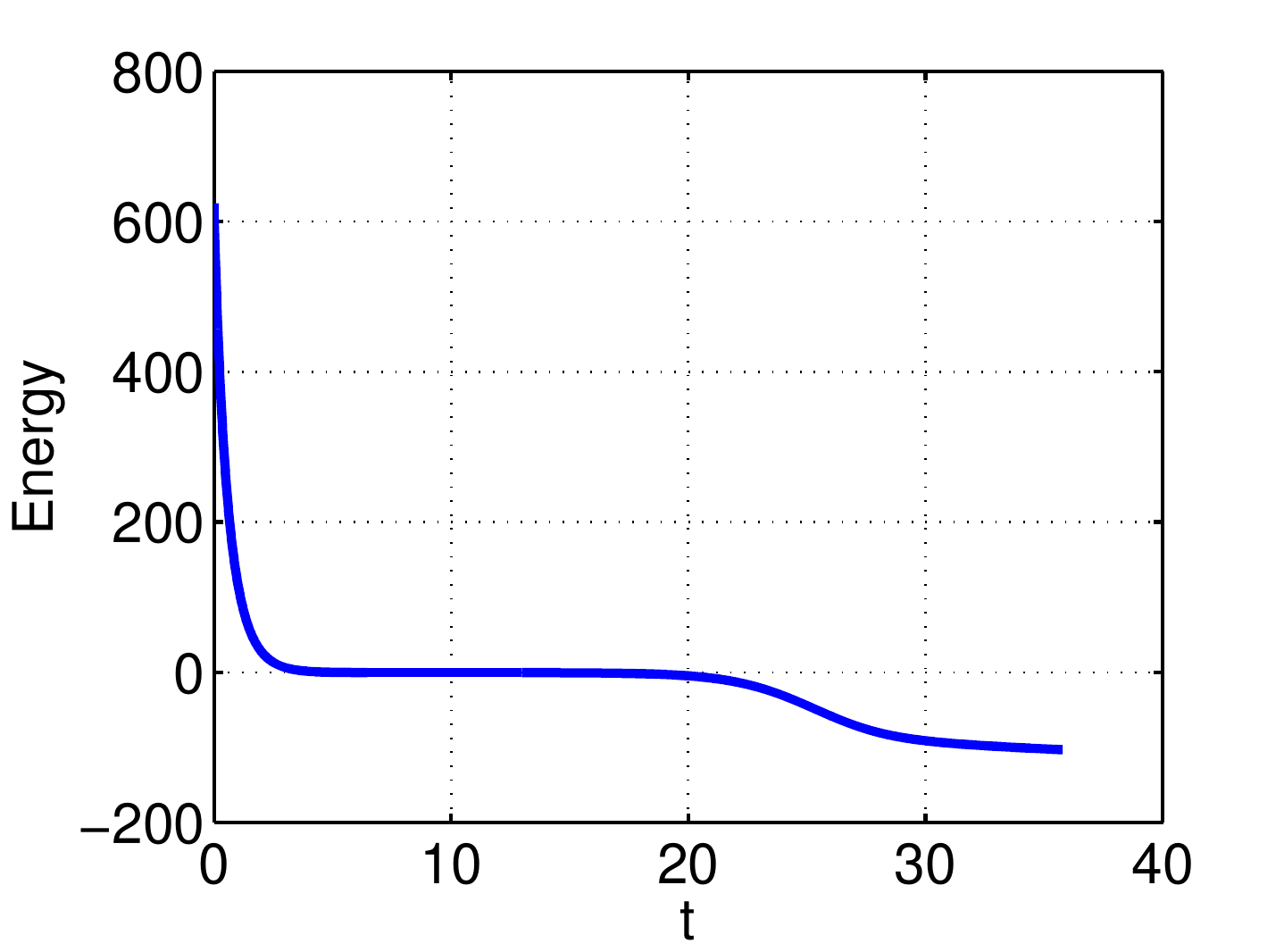}}
\subfloat{\includegraphics[scale=0.5]{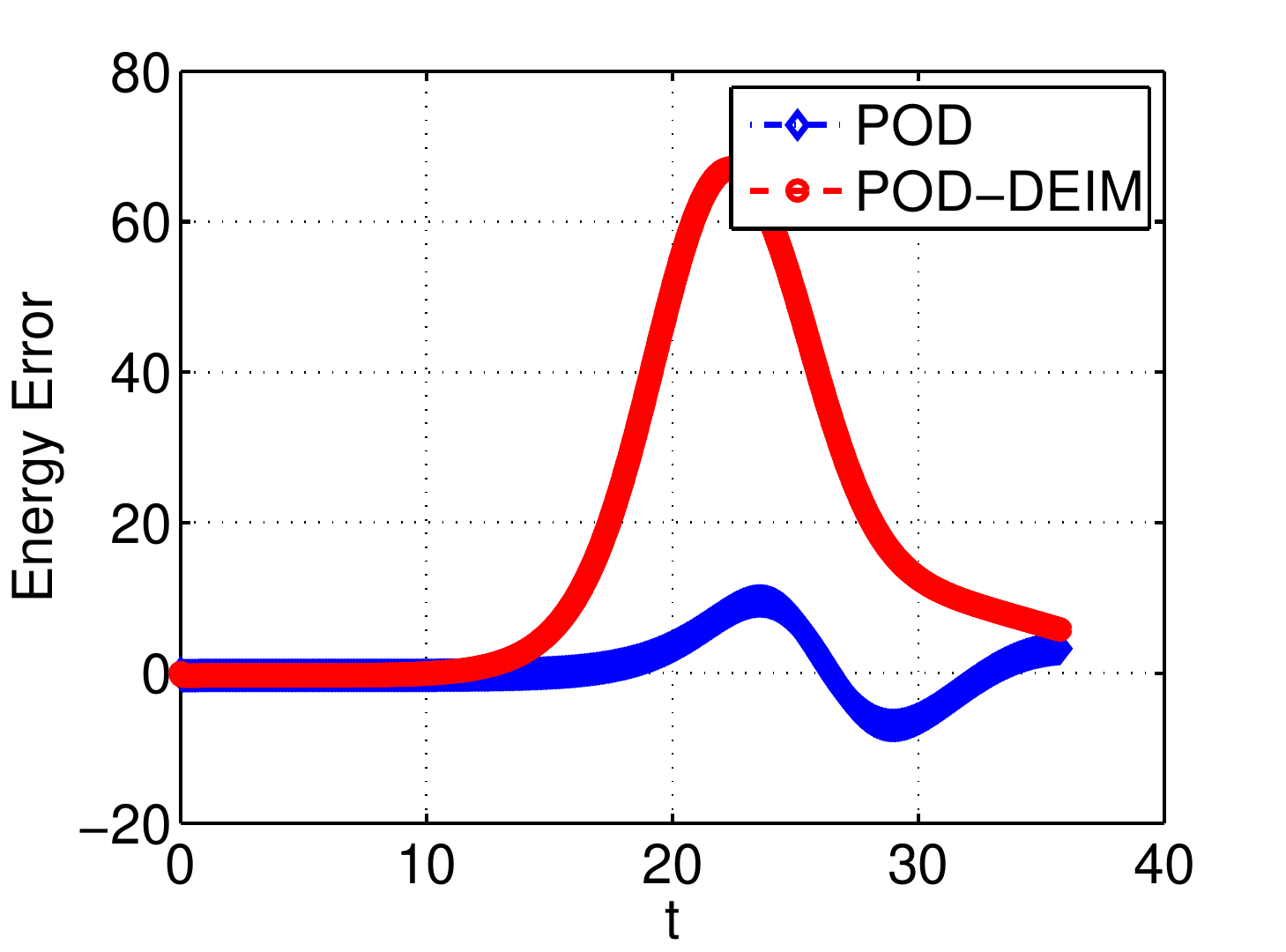}}
\caption{SH: Energy evolution of FOM  and energy errors of POD and DEIM\label{sh_energy}}
\end{figure}

Decay of the singular values shown in Fig.~\ref{sh_svd} is similar to the RGL equation in  Fig.~\ref{gle_svd}, and the energy dissipates like RGL equation as presented in  Fig.~\ref{sh_energy}. However, the energy decay is not well preserved for the POD and DEIM. Because the DEIM energy errors are larger then POD energy errors, the reduced order patterns obtained by DEIM are less accurate than those obtained by POD, see Table~\ref{sh_error} and Fig.~\ref{sh_rom}. Because more DEIM modes are required than the POD modes,  the speed-up of the DEIM is lower in Table \ref{sh_cpu} than for the RGL equation in  Table~\ref{gle_cpu}.

\begin{table}[H]
\centering
\caption{SH: POD \& DEIM errors\label{sh_error}}
\begin{tabular}{llllll}
\hline
    \# POD modes &     \#DEIM modes &      POD   &    DEIM  & POD & DEIM \\
  u(v)  & f & error &  error & energy error & energy error \\ \hline
       7(8)  &     15 &    8.73e-2      &  1.68e-01  &  1.09e-03  &  7.3e-03  \\
\hline
\end{tabular}
\end{table}

\begin{table}[H]
\centering
\caption{SH: CPU time and speedup factors\label{sh_cpu}}
\begin{tabular}{ l r r }
\hline
CPU time    & CPU time & speedup \\
\hline
FOM            & 779.32 & -  \\
POD  &  325.11  & 2.40 \\
DEIM  &   124.82   &  6.24\\
\hline
\end{tabular}
\end{table}

\section{Conclusions}

In this paper, we have developed a structure preserving reduced order modeling for the two typical gradient systems.
Numerical tests on two-dimensional real Ginzburg-Landau and Swift-Hohenberg equations  have illustrated the computational efficiency and the accuracy of the structure preserving reduced order models. The accuracy of the patterns and efficiency of the reduced order model depend on how well the energy evolution of the full order model is preserved  by  POD and DEIM. Because the singular values of the snapshot matrices and nonlinear terms decay slowly, relatively more POD and DEIM modes required.
The slow decay of the singular values is a characteristic for advection dominated and transport problems with wave type solutions. In the recent years, new methodologies are developed to tackle this problem.  Nevertheless, the time-dependent problems dominated by transport or propagation phenomena like the pattern formation remain a challenge for reduced order modeling. In a future study, we aim to apply these methodologies to reaction-diffusion systems with pattern formation.\\

\noindent {\bf Acknowledgments:}
The authors would like to thank the reviewer for the comments and suggestions that
helped to  improve the manuscript.


\end{document}